\documentclass[10pt, leqno]{article}
\usepackage{amsmath}
\usepackage{amsfonts}
\usepackage{amssymb}
\usepackage{graphicx}
\usepackage{color}

\newtheorem{thm}{Theorem}[section]
\newtheorem{cor}[thm]{Corollary}
\newtheorem{lem}[thm]{Lemma}
\newtheorem{prop}[thm]{Proposition}

\numberwithin{equation}{section}

\newcommand{\dx}{\,{\rm d}x}

\newcommand{\dt}{\,{\rm d}t}
\newcommand{\rd}{{\rm d}}

\newcommand{\dint}{\displaystyle\int}

\def\LL{\mathrm{L}} 
\newcommand{\RR}{\mathbb{R}}

\newcommand{\ve}{\varepsilon}
\def\ee{\mathrm{e}} 
\def\dist{\mathrm{dist}} 
\def\qed{\unskip\kern 6pt \penalty 500
\raise -2pt\hbox{\vrule \vbox to8pt{\hrule width 6pt
\vfill\hrule}\vrule}\par}
\begin{document}
\title{\bf Behaviour near extinction  for the Fast \\ Diffusion Equation on bounded domains}
\author{Matteo Bonforte$^{\,a}$,
Gabriele Grillo$^{\,b}$
~and~ Juan Luis Vazquez$^{\,c}$}

\maketitle

\begin{abstract}
 We consider the Fast Diffusion Equation $u_t=\Delta u^m$  posed in a bounded
smooth domain $\Omega\subset \RR^d$ with homogeneous Dirichlet
conditions; the exponent range is $m_s=(d-2)_+/(d+2)<m<1$. It is
known that bounded positive solutions $u(t,x)$ of such
problem extinguish in a finite time $T$, and also that such
solutions approach a separate variable solution $u(t,x)\sim
(T-t)^{1/(1-m)}S(x)$, as $t\to T^-$. Here we are interested in describing the  behaviour of the
solutions near the extinction time. We first show that the convergence
$u(t,x)\,(T-t)^{-1/(1-m)}$ to $S(x)$ takes place uniformly in the relative error norm. Then, we study the question of rates of convergence of the rescaled flow. For $m$ close to 1 we get such rates by means of entropy methods and weighted Poincar\'e inequalities. The analysis of the latter point makes an essential use of fine properties of the associated stationary elliptic problem $-\Delta S^m= {\bf c} S$ in the limit $m\to 1$, and such a study has an independent interest.
\end{abstract}

\vskip 1cm


\vfill

\noindent (a) Departamento de Matem\'{a}ticas, Universidad
Aut\'{o}noma de Madrid, Campus de Cantoblanco, 28049 Madrid, Spain.
E-mail address:{\tt~matteo.bonforte@uam.es}. \\
Web-page:{\tt~http://www.uam.es/matteo.bonforte}

\noindent (b) Dipartimento di Matematica, Politecnico di Milano, Piazza Leonardo da Vinci 32, 20183 Milano, Italy.
E-mail address:{\tt ~gabriele.grillo@polimi.it}

\noindent (c) Departamento de Matem\'{a}ticas, Universidad
Aut\'{o}noma de Madrid, and ICMAT, Campus de Cantoblanco, 28049 Madrid, Spain.
E-mail address:{\tt~juanluis.vazquez@uam.es}.\\ Web-page:{\tt~http://www.uam.es/juanluis.vazquez}

\newpage

\tableofcontents

\newpage

\section{Introduction}

We are interested in describing the  behaviour of nonnegative
solutions of the Fast Diffusion Equation (FDE) near the extinction
time. More precisely, we consider the following initial and boundary value
problem
\begin{equation}\label{FDE.Problem.Domain.Fine}
\begin{split}
\left\{\begin{array}{lll}
u_\tau=\Delta (u^m) & ~ {\rm in}~ (0,+\infty)\times\Omega\\
u(0,x)=u_0(x) & ~{\rm in}~ \Omega \\
u(\tau,x)=0 & ~{\rm for}~  \tau >0 ~{\rm and}~ x\in\partial\Omega\,,\\
\end{array}\right.
\end{split}
\end{equation}
posed in a bounded connected domain $\Omega\subset\RR^d$  with a regular boundary of class $C^{2,\alpha}$, $\alpha>0$. The fast diffusion range is $0<m<1$, but the theory developed below needs the further restriction $m_s<m<1$, where the lower end is the exponent  $m_s=(d-2)_+/(d+2)$ (inverse Sobolev exponent).
We assume that the initial data $u_0$ is a bounded and nonnegative function.
It is well-known that the above problem possesses a unique weak
solution $u\ge 0$ that is defined and positive for some time
interval $0<\tau<T$ and vanishes at $T=T(m,d,u_0)>0$, which is
called the (finite) extinction time, cf. \cite{BH, BV-ADV, DKV, VazLN}. Note that
the conditions on the initial data can be relaxed into
$u_0\in\LL^p(\Omega)$ for some $p\ge 1$ and $p>p_c$ where $p_c=\max\{1,
d(1-m)/2\}$ in view of the $\LL^p$-$\LL^\infty$ smoothing effect,
see \cite{BV-ADV, VazLN}.

\medskip

\noindent {\bf Rescaled equations and previous results.} To study the asymptotic behaviour it is convenient to transform the above problem by the known method of rescaling and time transformation. Thus, we put
\begin{equation}\label{transf}
u(\tau,x) = \left(\frac{T-\tau}{T}\right)^{\frac{1}{1-m}}v(t,x), \qquad
t= T\log\left(\frac{T}{T-\tau}\right).
\end{equation}
In this way, Problem {\rm (\ref{FDE.Problem.Domain.Fine})} is
mapped into the equivalent ``rescaled problem'':
\begin{equation}\label{FDE.Problem.Domain.Fine.Log}
\begin{split}
\left\{\begin{array}{lll}
v_{t}=\Delta (v^m)+\dfrac{v}{(1-m)T} & ~ {\rm in}~
(0,+\infty)\times\Omega,\\[3mm]
v(0,x)=u_0(x) & ~{\rm in}~ \Omega, \\[3mm]
v(t,x)=0 & ~{\rm for}~  t >0 ~{\rm and}~ x\in\partial\Omega.
\end{array}\right.
\end{split}
\end{equation}
The  transformation can also be expressed as
\begin{equation}\label{v}
v(t,x)
=\ee^{\frac{t}{1-m}}\, u\left(T-T\ee^{-t/T},x\right),
\end{equation}
and the time interval $0<\tau<T$ becomes $0<t<\infty$, so behaviour near extincion for the original flow becomes
behaviour as $t\to\infty$ in the rescaled flow, which is more convenient to analyze. Thus, in a
celebrated paper, Berryman and Holland \cite{BH}, 1980, reduced the study
of the behaviour near $T$ of the solutions of Problem
(\ref{FDE.Problem.Domain.Fine}) to the study of the possible
stabilization of the solutions of the transformed evolution
problem (\ref{FDE.Problem.Domain.Fine.Log}). They showed that the solutions
of the latter problem stabilize towards the solutions of the associated stationary
problem
\begin{equation}\label{FDE.Elliptic.Problem}
 -\Delta (S^m)=\frac{1}{(1-m)T}\,S \quad {\rm in}~ \Omega, \qquad S(x)=0 \quad {\rm for}~ x\in\partial\Omega\,,
\end{equation}
where $m_s<m<1$ and $\Omega$ are as before, and $S>0$ in $\Omega$. Using the new variable $V=S^m$ and putting
$p=1/m>1$ and ${\bf c}=1/((1-m)T)$ the latter problem can be written in the more popular semilinear elliptic form
\begin{equation}\label{eqn.V}
-\Delta V={\bf c}\,V^p \quad \mbox{in } \ \Omega, \qquad V=0 \quad \mbox{on } \ \partial\Omega
\end{equation}
Note that our restriction $m>m_s$ is the exact condition that makes the last problem subcritical, $p<p_s:=(d+2)/(d-2)$.

It is precisely proved in \cite{BH} that the rescaled orbit of a solution $v(t)=v(\cdot,t)$ converges in $W^{1,2}_0(\Omega)$ along subsequences to one or several stationary states $S$ (Remark: the elliptic problem can have multiple solutions depending on the geometry of $\Omega$). In the language of dynamical systems,  the omega-limit of $v$ is included in the set of positive classical solutions to the stationary problem \eqref{FDE.Elliptic.Problem}.

A main issue that remained open after \cite{BH} was to understand wether the rescaled solution $v$ converges to a unique stationary profile for all $t \to \infty$, even when the set of stationary solutions contains more than one function. The question of uniqueness of the asymptotic profile has been solved by Feireisl-Simondon in  \cite{FS2000}. We rewrite their main result in our notations:

\begin{thm}[\cite{FS2000}]\label{Thm.FS2000}
Let $\Omega\subset\RR^d$, $d\ge 1$ be a bounded domain of class $C^{2,\alpha}$, $\alpha>0$. Let $v\in\LL^{\infty}\big((0,\infty)\times\Omega\big)$ be a bounded weak solution of \eqref{FDE.Problem.Domain.Fine.Log}, then $v$ is continuous for $t>0$ and there exist a classical solution $S$ of the stationary problem \eqref{FDE.Elliptic.Problem} such that
\[
v(t)\to S \qquad\mbox{in } C(\overline{\Omega}) \qquad\mbox{as } t\to\infty\,.
\]
\end{thm}
We recall that Theorem 3.1 of \cite{FS2000} is a bit more general, indeed we specialize here to the case $f(u)=-u/(1-m)$.

We remark  that every solution $S=S_{m,T}(x)$ to the elliptic problem \eqref{FDE.Elliptic.Problem} produces a separable solution $\cal U$ of the original FDE
of the form
\begin{equation}\label{Separable.Solution}
\mathcal{U}(\tau,x)=S(x)\left(\frac{T-\tau}{T}\right)^{\frac{1}{1-m}}
\end{equation}
which corresponds to the initial datum
$\mathcal{U}(0,x)=S(x)$. Indeed, it is not a solution but a family of solutions since we can fix
$T>0$ at will: we will write ${\cal U}\,_T$ for definiteness when needed.

In the present paper we present  improvements on these results in several directions:

\noindent (i) We prove that the stabilization process takes place with Convergence in Relative Error. This topic occupies Section \ref{sec.cre}, and the main result is Theorem \ref{Thm.Main.Intro}.

\noindent (ii) In Sections \ref{sec.stab} and \ref{sec.Stabnear1} we prove convergence with rates  to the stationary state using so-called entropy methods. The use of the term entropy deserves an explanation: we introduce a suitable Lyapunov functional, namely a weighted $\LL^2$-norm of the quantity $v(t)-\overline{v}(t)$, that decreases in time along the nonlinear flow, and moreover, its dissipation in time is carefully controlled. We have decided to call entropy such functional since it has properties similar to the entropy functionals that have been extensively used to study the Cauchy problem, in particular the one used in \cite{BBDGV0, BBDGV, BDGV}; we do not claim a physical meaning for the entropy that we use here.

The first step is contained in Section \ref{sec.stab}, where we obtain convergence whenever a certain Weighted Poincar\'e Inequality holds with a sufficiently large constant, and precise decay is shown in that case, cf. Theorems \ref{asympt.rate1} and \ref{exp.decay.norm}.

\noindent(iii)  Then, in Section \ref{sec.Stabnear1} such assumption is  shown to hold for the solutions of our problem in  a restricted exponent range $m_\#<m<1$, and we obtain the concrete asymptotic results, Theorems \ref{thm.main.decay.rescaled} and \ref{rate.intro}.
The study relies heavily on the analysis of the associated semilinear equation \eqref{eqn.V} in the limit $p\to 1$, which has in our view an independent interest and is developed in a previous Section \ref{ell-sect}. The main result is Theorem \ref{Elliptic.THM}. In it we show that, if $V_p$ is a solution to equation \eqref{eqn.V} with homogeneous Dirichlet boundary datum  and we choose ${\bf c}=\lambda_p>0$ in such a way that $\|V_p\|_{p+1}=1$,  then $V_p/\Phi_1\to 1$ uniformly in $\Omega$ as $p\to 1$, where $\Phi_1$ is the ground state eigenfunction of the Dirichlet Laplacian on $\Omega$, with the normalization $\|\Phi_1\|_2=1$. Moreover we have that $\lambda_p\to \lambda_1$ as $p\to 1$ and, finally, Proposition \ref{prop.C} also shows that $1/(1-m)T\to \lambda_1$ as $m\to 1$.

\noindent (iv) The entropy method applies also for the Porous Medium Equation, that is when $m>1$, and it allows to find the rate of convergence, thus recovering the sharp result of Aronson and Peletier \cite{Ar-Pe}, by different methods. We devote Section \ref{Sect.PME} to present the the slow diffusion case $m>1$. It is worth mentioning that the we also obtain, as an intermediate result, a faster convergence for the entropy functional used in Section \ref{sec.stab}, see Theorem \ref{asympt.rate.PME} which is new.

\medskip

\noindent {\sc Notations.} Before proceeding with the statement and proofs of the results, let us recall some notations.   $S$ will denote the stationary solution of Problem \ref{FDE.Elliptic.Problem} indicated by Theorem \ref{Thm.FS2000}. For $x\in \Omega$ we write $d(x)=\dist(x,\partial\Omega)$ to indicate the distance to the boundary, its properties will de described below when needed.  $\lambda_1$ is the first eigenvalue of the Laplacian operator in the domain $\Omega$ with zero boundary conditions and positive eigenfunction $\Phi_1$. $\mathcal{S}_2$ is the optimal constant  in the Sobolev embedding $W^{1,2}_0(\Omega)\to L^{2^*}(\Omega)$, $2^*=2d/(d-2)$, $d\ge 3$. For $d\ge 3$ we put $p_s=2^*-1=(d+2)/(d-2)$
and $m_s=(d-2)/(d+2)=1/p_s$. By $\|\cdot\|_p$ we denote the standard $\LL^p(\Omega)$-norm, $1\le p\le \infty$, other norms will be carefully denoted.

\section{Convergence in relative error}\label{sec.cre}

We show that the quotient $v/S$ converges to 1 uniformly in the whole of $\Omega$, up to the boundary;

\begin{thm}\label{Thm.Main.Intro}
Let $u$ be the solution to Problem~{\rm \ref{FDE.Problem.Domain.Fine}} and let $T=T(m,d,u_0)$ be its extinction time. Let $S(x)$ be the positive classical solution to the  elliptic problem {\rm(\ref{FDE.Elliptic.Problem})}  predicted in Theorem~{\rm \ref{Thm.FS2000}}. Then,
\begin{equation}\label{eq.quot}
\lim_{\tau\to
T^-}\left\|\frac{u(\tau,\cdot)}{\mathcal{U}(\tau,\cdot)}-1\right\|_{\LL^{\infty}(\Omega)}
=\lim_{t\to\infty}\left\|\frac{v(t,\cdot)}{S(\cdot)}-1\right\|_{\LL^{\infty}(\Omega)}=0
\end{equation}
where $\mathcal{U}$ is the separable solution \eqref{Separable.Solution} and $v$ is  rescaled solution defined in \eqref{v}.
\end{thm}

This type of convergence is what we call {\sl uniform relative-error convergence} (REC for short),  and it is our first main contribution to the subject of fine asymptotics near extinction.  We can rephrase the result in terms of the rescaled solution $v$ as
\begin{equation}\label{LIMSUP.INF}
\lim_{t\to\infty}\left\|\frac{v^m(t,\cdot)}{S^m(t,\cdot)}-1 \right\|_{\LL^{\infty}(\Omega)}=0\,,
\end{equation}
and this form will be practical for the calculations. Since $V=S^m$  is a  function with linear growth near the boundary, what we say is that
\begin{equation}
\big[v^m(t,x) -S^m(x)\big]/d(x)\to 0
\end{equation}
uniformly in $x\in \Omega$ as $t\to\infty$.

As for related results,  DiBenedetto, Kwong and Vespri  proved in  \cite{DKV}, the Global Harnack Principle that we recall with our notations as follows, since we will be using it throughout this section:
\begin{thm}[\cite{DKV}]\label{GHP.DKV}
Let $v$ be the solution to the problem \eqref{FDE.Problem.Domain.Fine.Log} with $m_s<m<1$. Then, for any
$\varepsilon>0$ there exist positive constants $C_{0,m}\,,\,C_{1,m} >0$
depending on $d$, $m$, $\|u_0\|_{m+1}$, $\left\|\nabla
u_0^m\right\|_2$, $\partial\Omega$ and $\varepsilon$, such that for any
$t\ge \varepsilon$ and for any $x\in\Omega$
\begin{equation}\label{DKV.Estimates.Log}
C_{0,m}\, d(x)^{\frac{1}{m}}\le v(t,x) \le C_{1,m}\, d(x)^{\frac{1}{m}},
\end{equation}
Moreover for every $\kappa\ge 0$, $t\ge \varepsilon$ and for any $x\in\Omega$
\begin{equation}\label{DKV.analiticity}
\|v(t,\cdot)\|_{C^k(\Omega)}=\max_{|\alpha|=k}\,\sup_{x\in\Omega}\,|D^\alpha v(t,x)|
\le \frac{C_{1,m}^{|\alpha|+1} |\alpha|!}{ d(x)^{|\alpha|}}\,
d(x)^{\frac{1}{m}}\,.
\end{equation}
\end{thm}

Note that the constants $C_{0,m}\,,\,C_{1,m} >0$ may degenerate as $m\to 1$ or $m\to m_s$.
Estimate \eqref {DKV.Estimates.Log} above immediately implies the following estimate for the
solution to the FDE in the original variables:
\begin{equation}\label{DKV.Estimates}
C_{0,m}\, d(x)^{1/m}(T-\tau)^{1/(1-m)}\le u(\tau,x) \le C_{1,m}\,
d(x)^{1/m}(T-\tau)^{1/(1-m)}.
\end{equation}
This result is weaker than our Theorem \ref{Thm.Main.Intro} in the sense that our relative error convergence result is not a consequence of the above estimates, neither in the original nor in the rescaled variables: inequality \eqref{DKV.Estimates} only implies that the quotient $v(t)/S$ is bounded and bounded away from zero up to the boundary, but it does not prove that it converges to 1 as $t\to \infty$.
As far as we know, only the papers \cite{BH, DK, DKV, FS2000} contribute to the subject of the asymptotic of the Dirichlet problem for the FDE on bounded domains.

On the other hand, there is an extensive literature on stabilization of solutions of
evolution equations in different norms, mainly in $L^p$ or
$C^\alpha$ spaces. Let us comment on some
results on the topic of convergence in relative error which are
not so usual. Uniform convergence in relative error was first
proved for the Fast Diffusion Equation by one of the authors
\cite{Vascppme} in the following setting: solutions are
nonnegative and the equation is posed in the whole space $\RR^d$
with exponents $m_c<m<1$, where $m_c=(d-2)_+/d$. The case of all $m<1$ is treated in \cite{BBDGV, BGV, BDGV} where sharp rates of convergence in relative error are obtained for FDE posed in $\RR^d$ with a new entropy method, that does not apply to the bounded domain case.

In the present
setting of homogeneous Dirichlet data in a bounded domain, the
result is true for the Heat Equation (see a brief account in Section \ref{sec.lin}).
It is also true for the Porous Medium Equation, i.e.,
our problem for $m>1$, and the result follows from analyzing the
asymptotic result of Aronson and Peletier \cite{Ar-Pe}; see the
survey paper \cite{V2}. We recover the known results also in the PME case, which turns out to be simpler, see Section \ref{Sect.PME}\,. Hence, our present interest in the Fast
Diffusion case, where the presence of extinction makes the boundary argument more difficult, since the usual super-sub solution method does not work.

In the next subsections we proceed with the proof of Theorem~\ref{Thm.Main.Intro}.


\subsection{The relative error function and its equation}

\noindent From now on we will consider the evolution problem in its rescaled form
\eqref{FDE.Problem.Domain.Fine.Log}. Theorem \ref{Thm.FS2000} proves that the $\omega$-limit
of a solution is contained in the set of classical solutions to the Elliptic Problem
\eqref{FDE.Elliptic.Problem}, and the convergence takes place in
the uniform norm, and the solution $v$ selects a unique profile $S$ to converge to. Let us fix it: once we consider $v(0,t)=u_0$ then we know by Theorem \ref{Thm.FS2000} that $v(t)\to S$.

\noindent Now we introduce the {\it Relative Error Function} (REF)
\begin{equation}\label{def.phi.a}
\phi=\frac{v^m}{S^m}-1, \ \ \ \
v^m=S^m(\phi+1)=V(\phi+1),\qquad\mbox{and}\qquad V=S^m\,.
\end{equation}

\noindent $\bullet$~{\sl The Parabolic Equation of the REF and the
regularity of its solutions.}  Using the equations satisfied by $V$ and $v$, is
easy to show that $\phi$ satisfies the following parabolic
equation
\begin{equation}\label{Error.Parabolic.Equation.alpha}
\frac{1}{m}\left(1+\phi\right)^{\frac{1}{m}
-1}\phi_{t}=V^{1-\frac{1}{m}}\Delta\phi
+2\frac{\nabla V}{ V^{1/m}}\cdot\nabla\phi+F(\phi)
\end{equation}
where $F$ is given by
\begin{equation}\label{F.Error.Parabolic.Equation}
F(\phi)={\bf c}\left[\left(1+\phi\right)^{1/m}-\left(1+\phi\right)\right]
\end{equation}

\noindent Estimates \eqref{DKV.Estimates.Log} on $v$ and $S$, which is a stationary solution, imply that
\[
0<1-\left(\frac{C_{0,m}}{C_{1,m}}\right)^m=C_{2,m}\le
\phi\le C_{3,m}=\left(\frac{C_{1,m}}{C_{0,m}}\right)^m-1
\]
which proves that $\phi$ is bounded uniformly in
$(t,x)$, for $t>t_0>0$; notice that $t_0$ can be
chosen arbitrarily small, but this  affects the value of the
positive constants $C_{0,m}$ and $C_{1,m}$. Moreover, we notice
that
\[
1+\phi=v^m/ V>0
\]
in the interior of $\Omega$. Since  $\phi$ is also bounded
in the interior of $\Omega$, we conclude that the parabolic
equation {\rm (\ref{Error.Parabolic.Equation.alpha})} is neither
degenerate nor singular in the interior of $\Omega$. It follows
from standard quasilinear theory (cf. \cite{LSU}) that the
solution $\phi$ of such a parabolic equation is H\"older
continuous in any inner region $\overline{\Omega}_I\subset\Omega$
(the fact that $\phi$ is H\"older continuous could also be proved
by observing that $\phi+1=v^m/S^m$ and by recalling that both $v$
(see e.g. \cite{DKV}) and $S$ (see e.g. \cite{DKV} or
\cite{BH}) are at least H\"older continuous and positive in the
interior of $\Omega$).

\medskip

\noindent $\bullet$~{\sl Convergence of the REF in an interior
region of $\Omega$.} Under the running assumptions, we know by Theorem \ref{Thm.FS2000} that
\[
\sup_{\overline{\Omega}}\big|v(t)-S\big|\to 0\qquad\mbox{as}\qquad  t\to \infty\,,
\]
but this is not sufficient to prove the convergence
of the quotient $v^m/S^m$ to $1$ in the whole
$\Omega$, since at the boundary there is the problem caused by the
fact that both $v$ and $ V$ are zero and that the
parabolic equation \eqref{Error.Parabolic.Equation.alpha}
can degenerate at the boundary. However, such a problem is avoided
in any interior region where both $v$ and $S$ are strictly
positive. We define such interior region as
\[
\Omega_{I,\delta}=\left\{x\in\Omega\,:\,
\dist(x,\partial\Omega)>\delta\right\},
\]
with distance from the boundary $\delta>0$ which will be chosen later
small enough; we can thus say that in any interior region
$\Omega_{I,\delta}\subset\Omega$ we have
\begin{equation*}
\|\phi(t)\|_{\LL^{\infty}(\Omega_{I,\delta})}=
\sup_{\Omega_{I,\delta}}|\phi(t,\cdot)|\to
0,~ \mbox{~as~ }t\to\infty
\end{equation*}

recalling that $\phi=v^m/S^m -1$. We can sum up  what we have proved so far in the
following Lemma on Inner Convergence:

\begin{lem}\label{Lemma.Inner.Conv} Let $w$ be a solution to the
rescaled problem \eqref{FDE.Problem.Domain.Fine.Log}, and let $\phi$
be the corresponding relative error function defined by \eqref{def.phi.a}.
Then for every $\varepsilon>0$ and $\delta>0$ there exists
$t_{0,\varepsilon,\delta}>0$, such that for every
$t\ge t_{0,\varepsilon,\delta}$ and for every
$x\in\Omega_{I,\delta}$ we have
\begin{equation}\label{Inner.Convergence}
\left|\phi(t,x)\right|<\varepsilon\,.
\end{equation}
\end{lem}

\subsection{Distance to the boundary and barriers}

To get the proof of the convergence theorem we still have to show that uniform convergence of
$\phi$ takes place up to the boundary. To this end we will use a barrier argument, based on the following lemmas.
We remark here once and for all, that the barriers are independent of the particular choice of the stationary solution $S$.

\noindent First, we collect some properties of the function
``distance to the boundary''. It is defined as usual:
\[
d(x)=\dist(x,\partial\Omega)=\inf_{y\in\partial\Omega}|x-y|
\]
where $|\cdot|$ is the Euclidean norm of $\RR^d$.

\begin{lem}[Properties of the distance to the boundary]\label{Lemma.Distance}
Let $\Omega\subset\RR^d$ be a bounded domain with boundary
$\partial\Omega$ of class $C^2$. Let for $\xi>0$
\[
\Omega_{\xi}=\left\{x\in\Omega\,:\, d(x)<\xi
\right\}=\Omega\setminus \overline{\Omega_{I,\xi}}
\]
be the open strip of width $\xi$ near the boundary. Then,

\noindent {\rm (a)} there exists a constant $\xi_0>0$ such
that for every $x\in\Omega_{\xi_0}$, there is a unique
$h(x)\in\partial\Omega$ which realizes the distance:
\[
d(x)=|x-h(x)|.
\]
Moreover, $d(x)\in C^2(\Omega_{\xi_0})$ and for all $r\in[0,\xi_0)$ the function
$H_r:\partial(\overline{\Omega_r})\cap \Omega \to \partial\Omega$ defined by $H_r(x)=
h(x)$ is a homeomorphism.

\noindent {\rm (b)} Function $d(x)$ is Lipschitz with constant
$1$, i.e.
\[
|d(x)-d(y)|\le\, |x-y|.
\]
Moreover,
\[
0<c\le\,|\nabla d(x)|\le 1,~~~~\mbox{for any~}x\in\Omega_{\xi_0}
\]
and there exist a constant $K>0$ such that:
\begin{equation}\label{estimates.Laplacian}
-K\le\Delta d(x)\le K,~~~~\mbox{for any~}x\in\Omega_{\xi_0}
\end{equation}
\end{lem}

\noindent We refer to {\rm\cite{GT}} for the proof of this lemma.
Part (a) is due to Serrin.

We need a second technical result about estimates for the gradient
of the function $ V=S^m$ near the boundary. $S$ is a positive classical
solution to the elliptic problem {\rm
(\ref{FDE.Elliptic.Problem})}

 \begin{lem} \label{lem.grad.bry} For  $\xi_0>0$ small enough there exists
$\beta_0>0$ such that
\[
\nabla V(x)\cdot\nabla d(x) \ge \beta_0>0\,, \qquad \forall
x\in\Omega_{\xi_0}.
\]
\end{lem}

\noindent {\sl Proof.~} \noindent As explained in Section \ref{Section.Elliptic}, or as a consequence of the bounds \eqref{DKV.Estimates.Log}, the following estimate holds
\begin{equation}\label{estimates.sigma}
C_{0,m}^m d(x)\le  V(x)\le C_{1,m}^m d(x),~~~~\mbox{for
any~}x\in\Omega\,.
\end{equation}
Moreover, we know that $S$ is a positive classical solution of the elliptic
Dirichlet problem \ref{FDE.Elliptic.Problem} since
$m_s<m<1$. Thanks to Lemma \ref{Lemma.Distance} we can
conclude that both $ V$ and the distance function $d$ are functions
of class $C^2$ in a suitable neighborhood of the boundary
$\Omega_{\xi_0}$. The above estimates imply
\[
0\le C_{0,m}^m \frac{d(x)-d(x_0)}{|x-x_0|}\le
\frac{ V(x)- V(x_0)}{|x-x_0|}\le C_{1,m}^m
\frac{d(x)-d(x_0)}{|x-x_0|},
\]
for any $x\in\Omega_\xi$ and any $x_0\in\partial\Omega$, since
$d(x_0)= V(x_0)=0$, if $x_0\in\partial\Omega$. This implies
that
\[
0\le C_{0,m}^m \partial_j d(x_0)\le \partial_j  V(x_0) \le
C_{1,m}^m
\partial_j d(x_0)
\]
and that $\partial_j  V(x)$ and $\partial_j d(x)$ are both
nonnegative, so that $\nabla V(x)\cdot\nabla d(x) \ge 0$. Moreover
it implies that:
\[
C_{0,m}^m \sum_j(\partial_j d(x_0))^2\le \sum_j\partial_j
 V(x_0)\partial_j d(x_0) \le C_{1,m}^m \sum_j(\partial_j
d(x_0))^2
\]
and finally:
\begin{equation}\label{est.bound}
0<cC_{0,m}^m < C_{0,m}^m |\nabla d(x_0)|^2\le\nabla V(x_0)\cdot\nabla
d(x_0)\le C_{1,m}^m |\nabla d(x_0)|^2\le C_{1,m}^m
\end{equation}
for any $x_0\in\partial\Omega$, since we know from lemma {\rm (\ref{Lemma.Distance})} that
\[
0<c\le\,|\nabla d(x)|^2=\sum_j (\partial_j d(x))^2\le
1,~~~~\mbox{for any~}x\in\Omega_{\xi_0}\,.
\]
By continuity of $\nabla d(x)$ and since $|\Delta d(x)|\le K$, we can extend the estimates \eqref{est.bound} from $x_0\in\partial\Omega$, to a small neighborhood of the boundary, say $x\in\Omega_{\xi_0}$, eventually by putting a smaller lower constant $0<\beta_0\le cC_{0,m}^m $, which can eventually depend on $K>0$ and $\xi_0$\,.\qed

We have obtained uniform estimates on a neighborhood of the
boundary $\Omega_{\xi_0}$. These lemmas  are needed for a key ingredient in the proof of the
relative error convergence theorems, which is the
construction of barriers as super-solutions.

\begin{lem}\label{Lemma.Barriers.1}
We can choose positive  constants $A$, $B$, $C$ so that for every $t_0>0$ the function
\begin{equation}\label{barrier.Phi}
\Phi(t,x)=C-B\,d(x)-A(t-t_0)
\end{equation}
 is a super-solution to Equation
\eqref{Error.Parabolic.Equation.alpha} (the one satisfied by the REF $\phi$) on a
parabolic region
\[
\Sigma_{\Phi}=\Sigma_{\Phi,-1}=\left\{(t,x)\in(t_0,\infty)\times\Omega\,:\,\Phi(t,x)\ge
-1 \right\},
\]
and moreover $\Sigma_{\Phi}\subset
(t_0,\infty)\times\Omega_{\xi_1}$. Super-solution means that for
any $(t,x)\in \Sigma_{\Phi}$ we have
\[
\frac{1}{m}\left(1+\Phi\right)^{1/m -1}\Phi_{t}\ge
 V^{1-1/m}\Delta\Phi+2\frac{\nabla V}{ V^{1/m}}\cdot\nabla\Phi+F(\Phi)
\]
with
\[
F(\Phi)=\frac{1}{1-m}\left[\left(1+\Phi\right)^{1/m}-\left(1+\Phi\right)\right]\ge 0\\
\]
and $ V=S^m$, $S$ being a positive classical solution to the
elliptic problem {\rm (\ref{FDE.Elliptic.Problem})}. The constants
$A,B,C$ depend only on $m,d$, the upper bound for $ \xi_1$ and the
geometry of the border through $C_{0,m} $ and $C_{1,m}$.
\end{lem}
\noindent\textbf{Remark. } The construction of the barrier is quite technical, so we stress that when
$\xi_1$ is small enough, then a sufficient condition on the parameters is
\begin{equation}\label{cond1abc}
\left(1+C + \frac{1-m}m A\right)\,\xi_1\le (\beta B)^{m}.
\end{equation}
We considered the barrier on the parabolic region
$\Sigma_{\Phi,-1}$ since in that region the quantity $\Phi+1\ge
0$, but in what follows we will only  need the smaller region
$\Sigma_{\Phi,\varepsilon}\,$, for small $\varepsilon>0$.

\noindent {\sl Proof.~}We recall that the equation satisfied by the REF is
\begin{equation}\label{error.equation.lemma}
\frac{1}{m}\left(1+\phi\right)^{1/m
-1}\phi_{t}= V^{1-1/m}\Delta\phi+2\frac{\nabla V}{ V^{1/m}}
\cdot\nabla\phi+F(\phi).
\end{equation}
We will prove that function (\ref{barrier.Phi}) is
a super-solution for equation {\rm (\ref{error.equation.lemma})}
on the parabolic region $\Sigma_{\Phi}$ if we find constants $A$,
$B$ and $C$ such that
\begin{equation}\label{I.II}
\frac{1}{m}\left(1+\Phi\right)^{1/m -1}\Phi_{t}\ge (I)\ge
(II)\ge
 V^{1-1/m}\Delta\Phi+2\frac{\nabla V}{ V^{1/m}}\cdot\nabla\Phi+F(\Phi)
\end{equation}
where the expressions (I) and (II) are estimates from  below and
above respectively for the left and the right terms, independently
of $(t,x)\in\Sigma_{\Phi}$.

\noindent (I): This estimate is simple since  $\Phi_t  =-A$ for
any $(t,x)\in\Sigma_{\Phi}$. Hence:
\[
\frac{1}{m}\left(1+\Phi\right)^{1/m -1}\Phi_{t}=
-\frac{A}{m}\left(1+\Phi\right)^{1/m-1}=(I).
\]

\noindent(II): This estimate is more involved. First, we rewrite
the right-hand side of {\rm (\ref{I.II})} in a more convenient
form:
\[
\begin{split}
& V^{1-1/m}\Delta\Phi+2\frac{\nabla V}{ V^{1/m}}\cdot\nabla\Phi
+F(\Phi)=\frac{-B}{ V^{1/m}} \left[ V\,\Delta
d+2\nabla V\cdot\nabla d\right]+F(\Phi)\\&\le
\frac{-B}{ V^{1/m}} \left[ V\,\Delta
d+2\nabla V\cdot\nabla
d\right]+\frac{\left(1+\Phi\right)^{1/m}}{1-m}
\end{split}
\]
since
\[
\begin{split}
&\Delta\Phi = -B\Delta(d(x)),\qquad \nabla\Phi=-B\nabla(d(x)),\\
&F(\Phi)=\frac{1}{1-m}\left[(\Phi+1)^{1/m}-(\Phi+1)\right]\le\frac{1}{1-m}(\Phi+1)^{1/m}\\
\end{split}
\]
Moreover, we have that
\begin{equation}\label{lapl.grad.}
\left[ V\,\Delta d+2\nabla V\cdot\nabla d\right]\ge
-KC_{1,m}^m\xi_1+2\beta_0 =\beta>0
\end{equation}
on a region for any $x\in\Omega_{\xi_1}$ with
\begin{equation}\label{xi.lemma.barrier}
\xi_1=\min\left\{\xi_0,\frac{2\beta_0}{KC_{1,m}} \right\}
\end{equation}
with $K, \xi_0, \beta$ and $C_{1,m}>0$ as in the previous two
lemmas. Indeed, we have
\[
 V(x)\Delta d(x)\ge -K  V(x) \ge -KC_{1,m} \,d(x)\ge
-KC_{1,m}\xi_1
\]
for any $x\in\Omega_{\xi_1}$, as a consequence of estimate {\rm
(\ref{estimates.Laplacian})} and estimates {\rm
(\ref{estimates.sigma})}. Moreover,
\[
\nabla V(x)\cdot\nabla d(x) \ge \beta_0>0
\]
for any $x\in\Omega_{\xi_1}$, as proved in Lemma
\ref{lem.grad.bry}. Finally,
\[
\begin{split}
 V^{1-1/m} &\Delta\Phi +2\frac{\nabla V}{ V^{1/m}}
\cdot\nabla\Phi+F(\Phi)\le-\frac{B\,\beta}{\xi_1^{1/m}}+
\frac{\left(1+\Phi\right)^{1/m}}{1-m}=(II)\\
\end{split}
\]
in $\Omega_{\xi_1}$ with $\xi_1>0$ as in
{\rm(\ref{xi.lemma.barrier})}. With these estimates, we conclude
that $\Phi$ is a super-solution if the following condition holds:
\[
-\frac{A}{m}(1+\Phi)^{1/m-1}=(I)\ge(II)=-\frac{B\,\beta}{\xi_1^{1/m}}
+\frac{\left(1+\Phi\right)^{1/m-1}}{1-m}(1+\Phi).
\]
We can rewrite it in the form
\begin{equation}\label{Super.Condition.Tau}
\left[\frac{A}{m}+\frac{1+\Phi}{1-m}\right](1+\Phi)^{1/m-1}\le\frac{\beta\,B}{\xi_1^{1/m}}
\end{equation}
with $(t,x)\in\Sigma_{\Phi}$ and $\xi_1>0$ as in
{\rm(\ref{xi.lemma.barrier})}.

We have thus proved that $\Phi$ is a super-solution for any
$(t,x)\in\Sigma_{\Phi}$ with $x\in\Omega_{\xi_1}$ such that
{\rm (\ref{Super.Condition.Tau})} holds. Since $\Phi(t,x)\le C$
in the region under consideration, it suffices to choose $A$, $B$
and $C$ so that
\begin{equation}\label{Super.Condition}
\left[\frac{A}{m}+\frac{1+C}{1-m}\right](1+C)^{1/m-1}\le\frac{\beta\,B}{\xi^{1/m}}.
\end{equation}
to ensure that $(II)\le (I)$. \qed

\subsection{Proof of Theorem \ref{Thm.Main.Intro}}
\label{proof4.2}

It is based on the previous Lemmas and  Theorem \ref{Thm.FS2000} of \cite{FS2000}.

\medskip

\noindent (I) We have to show that  given $\varepsilon>0$ there
exists a time $T(\varepsilon)>0$ such that for any
$t>T(\varepsilon)$ and for any $x\in\Omega$ we have that
$|\phi(t,x)|<\varepsilon$. By Lemma \ref{Lemma.Inner.Conv} we
know that for every $\varepsilon>0$ and $\delta>0$ there exists
$t_{0,\varepsilon,\delta}>0$, such that for every
$t\ge t_{0,\varepsilon,\delta}$ and for every
$x$ in  interior region $\Omega_{I,\delta}$  we have
\[
\left|\phi(t,x)\right|<\varepsilon.
\]
It remains to show that uniform convergence takes place also up to
the boundary. This will be a consequence of comparison with the
barrier function of Lemma {\rm \ref{Lemma.Barriers.1}}, given by
\[
\Phi(t,x)=C-B\,d(x)-A(t-t_0).
\]
$A$, $B$, $C$ of the barrier $\Phi$ are suitable positive constants
which are chosen as in  Lemma \ref{Lemma.Barriers.1}, while
$t_0\ge 0$ is  a free parameter which will be adjusted later.

\medskip

\noindent (II)  Comparison  of $\phi$ with $\Phi$ takes place in a
neighborhood of the parabolic boundary of the form
$Q^*=(t_0,T)\times \Omega_\delta$ for $\delta$ small, and
$t_0,t_1$ to be determined. The parabolic border of this
region is formed  by three pieces: the initial section at
$t=t_0$, the inner parabolic boundary, and the outer lateral
boundary.  In order to compare $\phi$ and $\Phi$ we have to
check their values on the above three pieces of parabolic
boundary.

\noindent  \textsl{(a) We compare the values of $\phi$ and $\Phi$
at the initial section $t=t_0$.} We want that
\[
\phi(t_0,x)\le\Phi(t_0,x)=C-Bd(x)
\]
for all $x\in\Omega_{\delta}$. This is possible because of the
uniform boundedness of $\phi$
\[
C_{2,m}=c (C_{0,m}-C_{1,m})\le\phi(t_0,x)\le
C_{3,m}=c(C_{1,m}-C_{0,m})
\]
for all $x\in\Omega$ as a consequence of bounds
\eqref{DKV.Estimates.Log}. Now we simply choose $C$ sufficiently
large previous to the choice of $A$ and $B$ that have to satisfy
(\ref{cond1abc}).

\medskip

\noindent \textsl{(b) Comparison on the inner parabolic boundary:}
 This piece of the boundary is given by the points $(x,t)$ such that $d(x)=\delta$
and $t\in (t_0,T)$. On this set we want $\phi(t,x)\le
\Phi(t,x)$. Let us fix $\varepsilon>0$ and $0<\delta<\xi_1$
where $\xi_1>0$ is given in Lemma {\rm \ref{Lemma.Barriers.1}}. By
the uniform inner convergence (cf. Lemma \ref{Lemma.Inner.Conv})
we know that there exists a $t^*(\varepsilon,\delta)>0$ such
that $\phi(t,x)<\varepsilon$ on $\Omega_{I,\delta}$ if $t\ge
t^*$. The desired comparison holds if
\begin{equation}\label{cond2abc}
\varepsilon\le C-B\delta-A(t-t_0).
\end{equation}
Since $C$ cannot be small this implies restriction on  $B$ that
has to be compatible with (\ref{cond1abc}). This happens if
$\delta$ is small enough. Once $B$ and $C$ are chosen, it suffices
to take $A(t-t_0)$ small.

\medskip

\begin{figure}[ht]
\centering
\includegraphics[height=8cm, width=13cm]{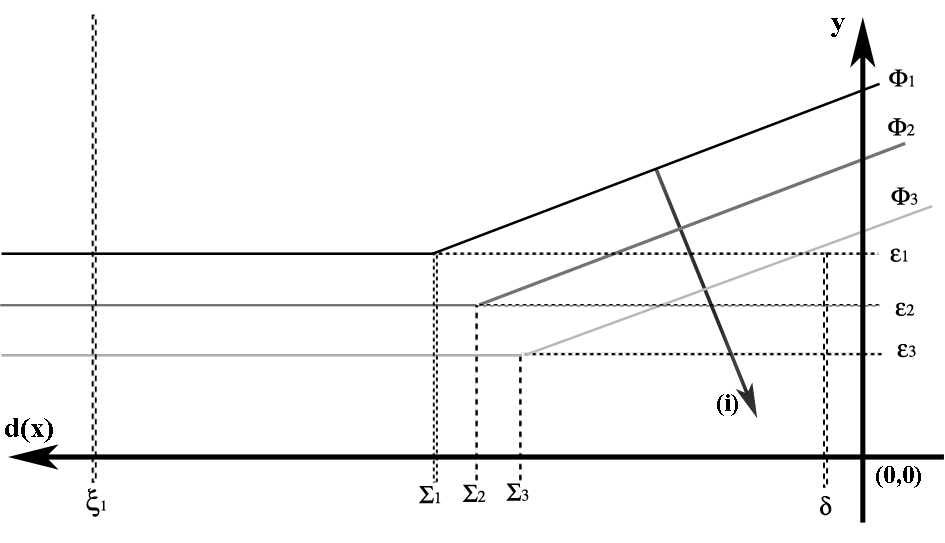}
\flushleft \noindent\textit{Figure 1\,: Idea of the behaviour of the barriers:\\
y-axis\,: values of $\Phi(t,x)$\\
x-axis\,: values of $d(x)=d(x,\partial\Omega)$, i.e. the distance from the boundary.\\
$\Sigma_i$\,: the points where $\Phi(t,x)=\varepsilon_i$,
i.e. the points of the boundary $\partial\Sigma_{\Phi,\varepsilon_i}\;$, $i=1,2,3.$\\
$\varepsilon_i$\,: different values of $\varepsilon$(decreasing
with $i=1,2,3$) give different barriers $\Phi_i$,\\ ~~~~decreasing with
$\varepsilon$ as the arrow (i) indicates.\\
$\xi_1$ and $\delta$ are as in Lemmas {\rm \ref{Lemma.Barriers.1}}
and {\rm \ref{Lemma.Barriers.3}.}}
\end{figure}

\medskip

\noindent \textsl{(c) We still have to check the comparison at the
outer lateral boundary}, $[t_0, T]\times \partial \Omega$,
where we only know that $\phi=v^m/S^m-1$ is bounded. But we can
use an approximation trick using the solutions $u_\eta$ of
problems posed in the domain $\Omega_\eta$ which is smaller than
$\Omega$. We know that $u_\eta\nearrow u$ as $\eta\to 0$.

Using $u_\eta$ instead of $u$ and $\Omega_\eta$ instead of
$\Omega$ allows to say that $\phi_\eta=u_\eta/\mathcal{U}-1=0$ on
the new outer boundary, hence $\phi_\eta<0$. Then we have:
\[
\phi_\eta=0<B\delta +\varepsilon < C-A(t-t^*)=\Phi
\]
thus $\phi_\eta<\Phi$ also on the outer boundary
$\partial\Omega_\eta$.

\medskip

Parabolic comparison allows then to say  that $\phi_\eta \le \Phi$
in the region $Q^*$ for $t\ge t_0\ge t^*$ such that
$t-t_0 <(C-B\delta-\varepsilon)/A$. Pass to the limit in
$\eta\to 0$ to get $\phi\le \Phi \quad \mbox{in} \ Q^*. $ In this
way  the following improvement of convergence near the boundary
after some time delay given by
\[
h_{\varepsilon,\delta}=\frac{C-B\delta-\varepsilon}{A},
\]
which is the maximum that \eqref{cond2abc} allows.

Steps (I) and (II) can be summarized in the following

\begin{lem}\label{Lemma.Barriers.3}
 Under the above conditions we have for  $t=t_0+h_{\varepsilon,\delta}$
\begin{equation}\label{Barrier.Bounds.Lemma.3}
\phi(t,x)\le\left\{\begin{array}{lll}
\varepsilon\,, &\mbox{for any ~}x\in\Omega,~\mbox{such that } \
d(x,\partial\Omega)>\delta,\\
\varepsilon+B\delta\,,  &\mbox{for any ~}x\in\Omega,~\mbox{such
that } \ d(x,\partial\Omega)<\delta.
\end{array}\right.
\end{equation}
 provided that $t_0\ge t^*$.
\end{lem}

\noindent (III) The proof of Theorem
\ref{Thm.Main.Intro} in the version of formula \eqref{LIMSUP.INF}  follows now by fixing
$\varepsilon>0$, finding a barrier with constants $A,B,C$ and then
taking $\delta < \varepsilon /B$. If $t\ge
t^*(\varepsilon,\delta)+ h_{\varepsilon,\delta}$, then
$$
\phi(t,x)\le2\varepsilon
$$
everywhere in $\Omega$. \qed

\subsection{The Porous Medium case}\label{Sect.PME}
\medskip

In this section we consider for the sake of comparison the same Dirichlet problem for the Porous Medium Equation
\[
\begin{split}
\left\{\begin{array}{lll}
u_\tau=\Delta (u^m) & ~ {\rm in}~ (0,+\infty)\times\Omega\\
u(0,x)=u_0(x) & ~{\rm in}~ \Omega \\
u(\tau,x)=0 & ~{\rm for}~  \tau >0 ~{\rm and}~ x\in\partial\Omega\,,\\
\end{array}\right.
\end{split}
\]
now for $m>1$. In this case there is no extinction in finite time, which changes the previous analysis and makes it much simpler. Let us give some details.  By means of the rescaling
\[
v(t,x)=u(\tau,x)(1+\tau)^{\frac{1}{m-1}}
\qquad\mbox{and}\qquad
1+\tau=\ee^t\,,
\]
the problem is mapped into the equivalent ``rescaled problem'':
\[
\begin{split}
\left\{\begin{array}{lll}
v_{t}=\Delta (v^m)+\dfrac{v}{m-1} & ~ {\rm in}~
(0,+\infty)\times\Omega,\\[3mm]
v(0,x)=u_0(x) & ~{\rm in}~ \Omega, \\[3mm]
v(t,x)=0 & ~{\rm for}~  t >0 ~{\rm and}~ x\in\partial\Omega.
\end{array}\right.
\end{split}
\]
The  transformation can also be expressed as
\[
v(t,x)
=\ee^{\frac{t}{m-1}}\, u\left(\ee^t-1,x\right) \qquad\mbox{with}\qquad t=\log(1+\tau)
\]
and the time interval $0<\tau<\infty$ remains $0<t<\infty$, in particular we preserve the initial datum. Notice that the rescaled problem when $m>1$ is formally the same as for the Fast Diffusion case, if one just considers ${\bf c}=1/(m-1)>0$. In this context it has been proved that solutions converge to a stationary state $S$, which is the unique solution of the elliptic equation $-\Delta S^m={\bf c}S$, see for example \cite{Ar-Pe, V2}. Solutions by separation of variables are given by $\mathcal{U}(\tau,x)=S(x)(1+\tau)^{-1/(m-1)}$\,.
The optimal rate of convergence in relative error as $\tau\to\infty$
\[
\left\|\frac{u(\tau,\cdot)}{\mathcal{U}(\tau,\cdot)}-1\right\|_{\LL^\infty(\Omega)}=\left\|\frac{u(\tau,\cdot)(1+\tau)^{\frac{1}{m-1}}}{S(\cdot)}-1\right\|_{\LL^\infty(\Omega)}\le \frac{K}{1+\tau}
\]
has been first obtained by Aronson and Peletier in \cite{Ar-Pe}, Theorem 3, for smooth nonnegative initial data. The rate can be easily shown to be optimal because of the special family of global solutions
\[
\mathcal{U}_{k}(\tau,x)=\frac{S(x)}{(k+\tau)^{\frac{1}{m-1}}}\qquad\mbox{for any }k>0\,.
\]
Indeed, the proof of \cite{Ar-Pe} strongly uses comparison with the stationary state $S$ through the special solutions $\mathcal{U}_k$. In the fast diffusion case comparison with smaller solutions does not help, since they extinguish earlier. This is just a example of the extra difficulties in proving convergence in relative error when $m<1$. In rescaled variables the result of \cite{Ar-Pe} reads
\begin{equation}\label{AP.exp}
\big| v(t,x)-S(x)\big|\le K\,S(x)\,\ee^{-t}\qquad\mbox{for all }x\in\Omega\;\mbox{and}\;t\gg 1\,.
\end{equation}
Other interesting approaches to convergence results, together with extensions to a larger class of initial data and solutions, can be found in \cite{V2}. Our entropy method applies to the case $m>1$ as we shall briefly discuss in the next sections, and allow us to recover the result of \cite{V2} and the optimal rates of convergence in relative error of \cite{Ar-Pe}, cf. Theorem \ref{PME.Rates}.

\section{Stabilization with convergence rates}\label{sec.stab}

 We recall here the setup and the notations. We consider the rescaled equation
\begin{equation}\label{FDE.Problem.Domain}\begin{split}
v_t = \Delta (v^m) &+ {\bf c}\,v,\\
{\bf c}=\frac{1}{T(1-m)}>0\quad\mbox{if }m<1, \qquad &\mbox{and}\qquad
{\bf c}=\frac{1}{m-1}>0\quad\mbox{if }m>1
\end{split}
\end{equation}
posed in a bounded connected domain $\Omega\subset\RR^d$ with sufficiently smooth boundary. We mainly deal with the so-called fast diffusion exponents and then we assume $m_s<m<1$. Almost all the calculations will hold for any $m>m_s$, including the case $m>1$; we will emphasize the differences when they will occur. The linear case $m=1$ is well known and will be briefly recalled in the next subsection as a motivation of our techniques.

We now introduce the quotient $w=v/S$, which converges uniformly to 1 on $\overline{\Omega}$, as a consequence of Theorem \ref{Thm.Main.Intro}. We then have
\begin{equation}
S w_t= \Delta (S^m w^m)+{\bf c} w\,,
\end{equation}
where $S$ is a stationary solution so that $\Delta S^m +{\bf c} S=0$ in $\Omega$ with $S=0$ on $\partial \Omega$, as precisely indicated in Theorem~\ref{Thm.FS2000}. We will also write $V =S^m$, which satisfies $\Delta V +{\bf c} V ^p=0$ in $\Omega$, $V =0$ on $\partial \Omega$ with exponent $p=1/m>1$.

We propose to perform the calculation on the asymptotic decay  in terms of $\theta=w-1=v/S-1$,  which is the relative error of the solution. Notice that it is different from the relative error $\phi$ used in the previous section, formula \eqref{def.phi.a}, which was defined as $\phi=w^m-1$, hence $\phi=(1+\theta)^m-1$.  We have
\begin{equation}
\theta_t= \frac1{S} \Delta (S^m (1+\theta)^m)+  {\bf c}\,(1+\theta).
\end{equation}
Using the identity
\begin{equation*}
\Delta (S^m (1+\theta)^m))= \nabla\cdot \left[ S^m \nabla (1+\theta)^m\right]+ \nabla (S^m)\cdot \nabla (1+\theta)^m
+ \Delta S^m \,(1+\theta)^m
\end{equation*}
and the equation for $S$, this  can be further written as
\begin{equation*}
\theta_t =\dfrac1{S} \nabla\cdot (S^m \nabla (1+\theta)^m) + \dfrac{\nabla S^m}{S}\cdot \nabla (1+\theta )^m
+{\bf c} \,f(\theta)\,,
\end{equation*}
where
\begin{equation}\label{f.taylor}
f(\theta):=(1+\theta)-(1+\theta)^m=(1-m)\left[\theta+\frac{m}{2}\theta^2+O\big(\theta^3\big)\right]
\end{equation}
for small $\theta$. One more calculation gives the following form for the equation
\begin{equation}\label{Eq.Theta}
S^{m+1}\theta_\tau
= \nabla\cdot (S^{2m }\nabla (1+\theta)^m) + {\bf c}\,S^{m+1} \,f(\theta)
\end{equation}
that we will use below. In terms of $V=S^m$ and $p=1/m$ we have
\begin{equation}
V ^{p+1}\theta_t
= \nabla\cdot (V ^{2 }\nabla (1+\theta)^m) + {\bf c}\,V ^{p+1} \,f(\theta)\,.
\end{equation}

\subsection{Weighted inequalities for the linear heat flow}\label{sec.lin}

In search for inspiration on how to proceed further, we compare the situation  with the standard way of treating linear equation $u_\tau=\Delta u$, which has striking formal similarities even if there is no extinction in finite time. After the suitable linear rescaling, which takes in the linear case the form $v(x,t)= e^{\lambda_1 t}v(x,t)$, we arrive at the same equation \eqref{FDE.Problem.Domain} with $m=1$ and ${\bf c}=\lambda_1$. The role of the stationary solution $S$ is now played by the first eigenfuction $\Phi_1>0$ of the Laplacian with zero boundary conditions (or any of its nonnegative multiples). The equation for $\theta=v/\Phi_1-1$ is:
\begin{equation}\label{eq.lin.theta}
\theta_t =\Delta \theta + 2\nabla\Phi_1\cdot\nabla \theta
        =\Phi_1^{-2}\nabla\cdot\left(\Phi_1^2\nabla \theta\right), \qquad \mbox{\rm i.\,e.,} \quad
        \Phi_1^2\,\theta_t=\nabla\cdot(\Phi _1^2\nabla \theta),
\end{equation}
to be compared with \eqref{Eq.Theta}. We recall next the known way to proceed with the asymptotic analysis of this equation via an Intrinsic Poincar\'e inequality. Indeed, let us consider the self--adjoint operator $-\Delta$ with Dirichlet boundary conditions on $\partial\Omega$. This can be defined by a general procedure by defining it to be the unique self--adjoint operator associated with the closure of the quadratic form $\int_\Omega|\nabla u|^2{\rm d}x/2$, initially defined for $u\in C_c^\infty(\Omega)$. Such operator has purely discrete spectrum, and we denote by $\lambda_j>0$, $j=1,2,\ldots$ its eigenvalues, arranged in nondecreasing order, and by $\Phi_j$ the corresponding $L^2$-normalized eigenfunctions.

The spectral representation for the corresponding heat semigroup $u_\tau=\Delta u$ shows that, letting $u_0$ be the initial datum and $c_j=\int_\Omega u_0\Phi_j\,{\rm d}x$, one has
$$
u(x,t)=\sum_{j=1}^{\infty}c_j e^{-\lambda_j t}\Phi_j(x)
$$
so that
\[
\theta:=\frac{u}{c_1\ee^{-\lambda_1 t}\Phi_1}-1\underset{t\to+\infty}\sim\frac{c_2}{c_1}\frac{\Phi_2}{\Phi_1}\ee^{-(\lambda_2-\lambda_1) t}.
\]
In other words, the solution $u(t)$ is close to the explicit solution $U_1(x,t)=c_1\ee^{-\lambda_1 t}\Phi_1$ and the relative error $w$ between such solutions, defined above, decays exponentially in time with a rate $\lambda_2-\lambda_1$. Notice in addition that the spatial factor ${\Phi_2}/{\Phi_1}$ is bounded.

To prepare the way to recovering a result of this kind in the nonlinear setting, where no spectral representation is available, we reformulate the above property as follows. Starting from equation \eqref{eq.lin.theta}, it is then natural to investigate the behaviour of $\theta$ by working in the weighted space L$^2(\Phi_1^2\,{\rm d}x)$. In fact first we observe that the weighted mean is preserved:
\[
\frac{\rd}{\dt}\int_\Omega \theta\Phi_1^2 \dx = \int_\Omega \nabla\cdot\left(\Phi_1^2\nabla \theta\right) \dx=0.
\]
Then we notice that:
\[
\frac{\rd}{\dt}\int_\Omega \theta^2 \Phi_1^2 \dx = 2 \int_\Omega \theta \nabla\cdot\left(\Phi_1^2\nabla \theta\right) \dx
    =-2\int_\Omega \left|\nabla \theta\right|^2\Phi_1^2\dx\,.
\]
By the the above conservation of weighted mean we can and shall assume that ${\theta}_{\Phi_1}=0$, where
\[
{g}_{\Phi_1}=\frac{\int_\Omega g \Phi_1^2 \dx}{\int_\Omega\Phi_1^2 \dx}\,.
\]
It is then clear that in order to get a decay rate for $E[\theta]=\int_\Omega \theta^2 \Phi_1^2 \dx$ it suffices to prove the following \it intrinsic Poincar\'e inequality:\rm
\begin{prop}\label{Intrinsic.Poincare}
Let $f\in W_0^{1,2}(\Omega)$ and  $g=f/\Phi_1$. Then the following inequality
holds
\begin{equation}\label{Poincare.w}
(\lambda_2-\lambda_1)\int_\Omega \left|g-{g}_{\Phi_1}\right|^2 \Phi_1^2 \dx
    \le \int_\Omega \left|\nabla g\right|^2\Phi_1^2\dx.
\end{equation}
\end{prop}

Although this inequality is well-known, we provide in Appendix \ref{sec.appendix}.2  a short proof for the reader's convenience, and we recall there some sharp upper and lower bounds on the spectral gap $\lambda_2-\lambda_1$.

\subsection{Energy analysis of the nonlinear flow}

Inspired by the preceding linear analysis  and after carefully choosing among the different options to
attack the nonlinearities of our evolution process, we are going to prove  a certain type of entropy/entropy-production inequalities. We define the suitable entropy functional to be
\begin{equation}\label{E.theta}
\mathcal{E}[\theta(t)]=\frac{1}{2}\int_\Omega\big|\theta(t)-\overline{\theta}(t)\big|^2 S^{1+m}_{c,m}\dx\,,
\end{equation}
where
\[
\overline{\theta}(t)=\frac{\int_\Omega\theta(t,x)S_{c,m}^{1+m}\dx}{\int_\Omega S^{1+m}_{c,m}\dx}\,.
\]
and  $S=S_{c,m}$ is the chosen solution to the elliptic problem
\begin{equation}\label{ellipt}
\left\{
\begin{array}{lll}
-\Delta S^m={\bf c}\,S&\mbox{in }\Omega\\
S>0&\mbox{in }\Omega\\
S=0&\mbox{on }\partial\Omega\,,
\end{array}
\right.
\end{equation}
 ${\bf c}>0$ and $m_s<m<1$, and the relative error $\theta$ satisfies the equation
\begin{equation}\label{eq.theta}
\theta_t
= \frac{1}{S^{m+1}} \nabla\cdot (S^{2m }\nabla (1+\theta)^m) + {\bf c} \,f(\theta),\qquad\mbox{with}\qquad f(\theta)=(1+\theta)-(1+\theta)^m\,.
\end{equation}
We do not specify boundary conditions, but we know that $\theta$ is continuous up to the boundary $\partial\Omega$ and the convergence in relative error valid in $C(\overline{\Omega})$, proved in Theorem \ref{Thm.Main.Intro} indicates that the boundary conditions stabilize to 0 as $t\to\infty$. It is not restrictive to assume that $|\theta|\le\varepsilon$ on $\partial\Omega$, or by the maximum principle in the whole $\overline{\Omega}$, for arbitrarily small $\varepsilon>0$. The price that we have to pay is just a time shift. We now prove the following

\begin{prop}\label{prop.deriv.entr}
Let ${\bf c}>0$ be as in \eqref{FDE.Problem.Domain}, and $m_s<m<1$. Let $\theta$ be a global smooth solution to equation \eqref{eq.theta} and let $\varepsilon(t):=\|\theta(t,\cdot)\|_\infty\to 0$. Then,  the following inequality holds
\begin{equation}\label{Entr.Prod.}\begin{split}
-\frac{\rd}{\rd t}\mathcal{E}[\theta(t)]
    &\ge m[1+\varepsilon(t)]^{m-1}\int_\Omega \,|\nabla \theta(t,x)|^2\,S^{2m}\dx
    - 2{\bf c}[1-m+\varepsilon(t)]\, \mathcal{E}[\theta(t)]\\
\end{split}
\end{equation}
for all times $t>t_0$ where $t_0$ is such that $\varepsilon(t)=\|\theta(t,\cdot)\|_\infty<1$ for all $t\ge t_0$. When $m=1$ we recover the standard $\LL^2$-weighted estimates that hold in the linear case, see Subsection \ref{sec.lin}. Moreover, when $m>1$, the property $\varepsilon(t)\to 0$ holds and, moreover, for $t$ sufficiently large:
\begin{equation}
-\frac{\rd}{\dt}\mathcal{E}[\theta(t)]
    \ge m[1-\varepsilon(t)]^{m-1}\int_\Omega \,|\nabla \theta|^2\,S^{2m}\dx
    + 2{\bf c}[m-1-\varepsilon(t)]\, \mathcal{E}[\theta(t)]\ge 0\,.
\end{equation}
\end{prop}
\noindent {\sl Proof.~} Uniform convergence in relative error holds as $t\to\infty$ by Theorem \ref{Thm.Main.Intro}, that means $\|\theta(t,\cdot)\|_\infty\to 0$. We set now $\rd\mu=S_{c,m}^{1+m}\dx$ and we will write $S=S_{c,m}$ throughout the proof, since no confusion will arise. First we notice that
\begin{equation}\label{mean}
\begin{split}
\int_\Omega\big[\theta(t)-\overline{\theta}(t)\big]\,\big[\partial_t\overline{\theta}(t)\big]\,\rd\mu
& =\partial_t\overline{\theta}(t)\int_\Omega\big[\theta(t)-\overline{\theta}(t)\big]\,\rd\mu =0\,,
\end{split}
\end{equation}
and also that
\begin{equation}\label{var}
0\le \int_\Omega \big|\theta-\overline{\theta}(t)\big|^2 \rd\mu
    =\int_\Omega \big[\theta^2-\overline{\theta}(t)^2\big]\rd\mu\,.
\end{equation}
We next differentiate $\mathcal{E}[\theta(t)]$ along the flow
\[
\begin{split}
\frac{\rd}{\dt}\mathcal{E}[\theta(t)]
    &=\int_\Omega\big[\theta(t)-\overline{\theta}(t)\big]\,\big[\partial_t\theta(t)\big]\,\rd\mu
        +\int_\Omega\big[\theta(t)-\overline{\theta}(t)\big]\,\big[\partial_t\overline{\theta}(t)\big]\rd\mu\\
    &=-\int_\Omega \nabla\left[\theta(t)-\overline{\theta}(t)\right]
            \cdot\left[S^{2m} \nabla(1+\theta)^m\right]\dx
    +{\bf c}\int_\Omega \left[\theta(t)-\overline{\theta}(t)\right]\,f(\theta)\,\rd\mu\\
    &=-m\int_\Omega (1+\theta)^{m-1}\,|\nabla \theta|^2\,S^{2m}\dx
        +{\bf c}\int_\Omega \left[\theta(t)-\overline{\theta}(t)\right]\,f(\theta)\,\rd\mu=(I)+(II)\,.\\
\end{split}
\]
In order to estimate $(II)$ we notice that, since $\overline{\theta}$ does not depend on the spatial variable:
\[
\int_\Omega f\big(\overline{\theta}(t)\big)\big(\theta-\overline{\theta}(t)\big)\rd\mu
=f\big(\overline{\theta}(t)\big)\int_\Omega\big(\theta-\overline{\theta}(t)\big)\rd\mu=0\,.
\]
We also notice that  for small $\theta-\overline{\theta}$ we have
\[
f(\theta)=f(\overline{\theta})+f'(\overline{\theta})(\theta-\overline{\theta}) + \frac12 f''(\widetilde{\theta}) (\theta-\overline{\theta})^2
\]
where $\widetilde{\theta}$ lies between $\theta$ and $\overline{\theta}$. We then have
\[
\begin{split}
\frac{(II)}{{\bf c}}&=\int_\Omega f(\theta)\left[\theta-\overline{\theta}(t)\right]\,d\mu
        =\int_\Omega \left[f(\theta)-f\big(\overline{\theta}(t)\big)\right]\big(\theta-\overline{\theta}(t)\big)\rd\mu\\
    &=f'(\overline{\theta})\int_\Omega\left[\theta-\overline{\theta}(t)\right]^2\,d\mu
    + \frac12  \int_\Omega f''(\widetilde{\theta}) \left[\theta-\overline{\theta}(t)\right]^3\,d\mu.
\end{split}
\]
Now we use the fact that
\[
f'(\overline{\theta})=1-m(1+\overline{\theta})^{m-1}=1-m + m(1-m)\overline{\theta}+ O(\overline{\theta}^2)
\]
which tends to $(1-m)$ as $t\to \infty$ uniformly in the space. Also
$$
f''(\widetilde{\theta}) = m(1-m)(1+\widetilde{\theta})^{m-2}\to m(1-m)
$$
uniformly in space as $t\to\infty$. Putting these things together, we
have obtained that
\begin{equation}\label{Entr.Prod.1}\begin{split}
-\frac{\rd}{\dt}\mathcal{E}[\theta(t)]
    &\ge m\int_\Omega (1+\theta)^{m-1}\,|\nabla \theta|^2\,S^{2m}\dx
    - 2{\bf c}[1-m+\varepsilon(t)]\, \mathcal{E}[\theta(t)]\\
    &\ge m[1+\varepsilon(t)]^{m-1}\int_\Omega \,|\nabla \theta|^2\,S^{2m}\dx
    - 2{\bf c}[1-m+\varepsilon(t)]\, \mathcal{E}[\theta(t)]\\
\end{split}
\end{equation}
notice that in the limit $m\to 1$, the last term disappears, and we recover the standard $\LL^2-$weighted estimates that hold in the linear case, see subsection \ref{sec.lin}.  When $m>1$ we obtain, using convergence in relative error as proved in \cite{Ar-Pe}:
\begin{equation}\begin{split}
-\frac{\rd}{\dt}\mathcal{E}[\theta(t)]
    &\ge m\int_\Omega (1+\theta)^{m-1}\,|\nabla \theta|^2\,S^{2m}\dx
    - 2{\bf c}[1-m+\varepsilon(t)]\, \mathcal{E}[\theta(t)]\\
    &\ge m[1-\varepsilon(t)]^{m-1}\int_\Omega \,|\nabla \theta|^2\,S^{2m}\dx
    - 2{\bf c}[1-m+\varepsilon(t)]\, \mathcal{E}[\theta(t)]\\
\end{split}
\end{equation}
as claimed. The fact that the r.h.s. is positive for sufficiently large time follows again by convergence in relative error.\qed

\subsection{Weighted Poincar\'e Inequality and first rate of convergence}\label{ssec.frc}

In order to get a rate of decay for $\mathcal{E}[\theta(t)]$ we shall need a suitable version of the weighted Poincar\'e inequality adapted to our problem, that we formulate next:

\noindent {\bf GWPI: General Weighted Poincar\'e inequality.} \emph{Given $m,\Omega$,  there exists a constant $K>0$ such that  for every stationary solution $S>0$ of problem \eqref{ellipt} with constant ${\bf c}=1$, and for every $\theta\in W^{1,2}_0(\Omega,S^{2m}dx)$
we have}
\begin{equation}
 \int _\Omega S^{2m}|\nabla \theta|^2\,dx \ge K \int_\Omega S^{m+1}|\theta-\overline \theta|^2\,dx\,.
\end{equation}
Note that $K$ depends only on $m$ and $\Omega$. In order to apply this property to positive solutions $S_c=S_{m,c}$ of the elliptic problem \eqref{ellipt} with constant ${\bf c}\ne 1$
we use the transformation
\begin{equation}
S_c(x)=\mu \, S_1(x)
\end{equation}
that produces a solution $S_1$ of problem \eqref{ellipt} if $\mu={\bf c}^{1/(1-m)}$. Therefore, we have
\begin{equation*}
\begin{array}{l}
 \dint _\Omega S_c^{2m}|\nabla \theta|^2\,dx =\mu^{2m}\dint _\Omega S_1^{2m}|\nabla \theta|^2\,dx
 \ge K \mu^{2m} \dint_\Omega S_1^{1+m}|\theta-\overline \theta|^2\,dx=K\mu^{m-1}\dint_\Omega S_c^{m+1}|\theta-\overline \theta|^2\,dx\,.
 \end{array}
\end{equation*}
In conclusion, the GWPI is formulated for problem \eqref{ellipt} with constant ${\bf c}\ne 1$ as
\begin{equation}\label{ewpi2}
 \dint _\Omega S_c^{2m}|\nabla \theta|^2\,dx \ge K{\bf c} \dint_\Omega S_c^{m+1}|\theta-\overline \theta|^2\,dx\,.
\end{equation}

\noindent {\bf Main assumption.} We now make an assumption that is crucial for the rest of the paper
\begin{equation}\label{K-assump}
Km-2(1-m)\ge \lambda_0>0.
\end{equation}

\noindent The rest of the paper will be based on deriving the consequences of this assumption and on the other hand, on justifying that under suitable conditions the assumption holds. As a first hint in the latter direction, if we consider the linear case and compare formula \eqref{Poincare.w} in Appendix \ref{sec.appendix}.3, with formula  \eqref{ewpi2} we see that the latter holds with $m=1$,  ${\bf c}=\lambda_1$ and $K=(\lambda_2-\lambda_1)/\lambda_1$, hence $\lambda_0=K>0$ in \eqref{K-assump}, and finally $\lambda_0 {\bf c}=\lambda_2-\lambda_1$. When we deal with the PME case, this assumption is always satisfied since $m>1$, hence we do not care about the expression of the constant $K$ in the GWPI.

Here is a first important consequence of that assumption in the FDE case:

\begin{thm}\label{asympt.rate1} Assume  that \eqref{K-assump} holds for a given $m\in (m_s,1)$ and given $\Omega\subset\RR^d$. Let $\theta$ be a global bounded and positive solution to equation \eqref{eq.theta}. Then for every $\gamma<\gamma_0:=\lambda_0 {\bf c}$ there exists a time $t_0>0$ such that for all $t\ge t_0$
\begin{equation}
 \qquad \qquad \mathcal{E}[\theta(t)]\le \mathcal{E}[\theta(t_0)]e^{-\gamma  (t-t_0)}.
\end{equation}
\end{thm}

\noindent {\sc Proof.} It is immediate after the derivation so far.
Assuming now that \eqref{K-assump} holds we go back to formula  \eqref{Entr.Prod.1} to get
$$
-\frac{\rd}{\dt}\mathcal{E}[\theta(t)]
    \ge {\bf c} \left(K m[1+\varepsilon(t)]^{m-1}
    - 2[1-m+\varepsilon(t)]\right) \, \mathcal{E}[\theta(t)].
$$
As $\varepsilon(t)\to 0$ when $t\to \infty$, the conclusion holds. \qed

Next we deal with the simpler PME case where we just need a GWPI with any constant to get a rate.
\begin{thm}\label{asympt.rate.PME} Let $m>1$ and let $\theta$ be a global bounded and positive solution to equation \eqref{eq.theta}. Then, for all $\beta<2+\frac{Km}{m-1}$ there exists a time $t_1$ depending on $m,d,\beta$ and on the constant $K>0$ of the GWPI, such that
\begin{equation}
\mathcal{E}[\theta(t)]\le \mathcal{E}[\theta(t_1)]\,\ee^{-\beta(t-t_1)}\qquad\mbox{for all }t\ge t_1.
\end{equation}
\end{thm}
\noindent {\sc Proof.} It is immediate after the derivation so far. Recall that ${\bf c}=1/(m-1)>0$ in this case.  We go back to formula \eqref{Entr.Prod.} to get, for any $\delta>0$ and any $t\ge t_\delta$ large enough:
\[\begin{split}
-\frac{\rd}{\dt}\mathcal{E}[\theta(t)]
    &\ge  \frac{1}{m-1}\left(K m[1-\varepsilon(t)]^{m-1}
    + 2[m-1-\varepsilon(t)]\right) \, \mathcal{E}[\theta(t)] \\
    &=\frac{1}{m-1}\left[K m[1-\varepsilon(t)]^{m-1}-2\varepsilon(t)
    + 2[m-1]\right] \, \mathcal{E}[\theta(t)]\ge \left[2+\frac{Km}{m-1}-\delta\right]\, \mathcal{E}[\theta(t)]\mbox{.\qed}
\end{split}\]
\subsection{Norm decay}

There is one more step to perform, since the entropy decay does not automatically imply the decay of the weighted $\LL^2$-norm, because the mean value $\overline{\theta}$ is not constant along the nonlinear evolution under consideration hence it must be also controlled. This is also related to the fact that we can not specify boundary conditions for the relative error function $\theta$. We first deal with the main case $m_s<m<1$. Recall that $S=S_{c,m}.$
\begin{thm}\label{exp.decay.norm}
Under the assumptions of Theorem {\rm \ref{asympt.rate1}}. Then both entropy and the $\LL^2-$norm, decay exponentially with the same rate $\gamma<\gamma_0=\lambda_0 {\bf  c}$. More precisely there exists a constant $\kappa>0$ such that
\begin{equation}\label{rate.entropy}
\int_\Omega\big|\theta(t)\big|^2 S^{1+m} \dx\le \kappa\,\mathcal{E}[\theta(t)]\le \kappa\,\ee^{-\gamma(t-t_1)}\mathcal{E}[\theta(t_1)]
\end{equation}
for all $t>t_1\gg 1$, where $\kappa$ depends on $m,d$ and $\mathcal{E}[\theta(t_1)]$, cf. the end of the proof.
\end{thm}
\noindent {\sl Proof.~}We first observe that since $\|\theta(t)\|_\infty\to 0$ as $t\to \infty$, hence also $\overline{\theta}(t) \to 0$ as $t\to \infty$,  we can always assume $|\theta|<1/2$ and $|\overline{\theta}|<1/2$, for all $t\ge t_0$. Next, we deduce the differential equation for
\[
\overline{\theta}(t)=\frac{\int_\Omega\theta(t,x)S ^{1+m}\dx}{\int_\Omega S^{1+m} \dx}\,,
\]
using the equation $\theta_t
= S ^{-(m+1)}\nabla\cdot (S ^{2m }\nabla (1+\theta)^m) + {\bf c} \,f(\theta)$:
\begin{equation*}\begin{split}
\overline{\theta}\,'(t)
    &=\frac{1}{\int_\Omega S^{1+m} \dx}\frac{\rd\overline{\theta}}{\dt}(t)
     =\frac{1}{\int_\Omega S^{1+m} \dx}\frac{\rd}{\dt}\int_\Omega \theta(t,x)S ^{m+1}\dx\\
    &= \frac{1}{\int_\Omega S^{1+m} \dx}\int_\Omega \nabla\cdot (S ^{2m }\nabla (1+\theta)^m)\dx
        + \frac{{\bf c}}{\int_\Omega S^{1+m} \dx} \int_\Omega  f(\theta)S ^{m+1}\dx\\
    &= \frac{{\bf c}\,\int_\Omega  f(\theta)S ^{m+1}\dx}{\int_\Omega S^{1+m} \dx} \\
\end{split}
\end{equation*}
where $f(\theta)=(1+\theta)-(1+\theta)^m$. By convexity of $f$ we have that $f$ lies above its tangent at the origin,
\[
f(\overline{\theta})\ge (1-m)\overline{\theta}\,,
\]
so that
\[
\overline{\theta}\,'(t)= \frac{{\bf c}\,\int_\Omega  f(\theta)S ^{m+1}\dx}{\int_\Omega S^{1+m} \dx}\ge {\bf c} f(\overline{\theta}(t))\ge {\bf c}(1-m)\overline{\theta}(t)\,,
\]
where in the first step we have used Jensen's inequality since $f$ is convex. An integration over $[s,t]\subset [t_0,\infty)$ gives
\[
\overline{\theta}(t)\ge \overline{\theta}(s)\ee^{{\bf c}(1-m)(t-s)}
\]
which implies that $\theta(s)\le 0$ for all $s\ge t_0$, otherwise we get a contradiction with the fact that $\overline{\theta}(t)\to 0$ as $t\to +\infty$. Next by Taylor expansion we get
\begin{equation}\label{ineq.f}
f(\theta)=f(\overline{\theta})+f'(\overline{\theta})(\theta-\overline{\theta})
        + \frac12 f''(\tilde{\theta}) (\theta-\overline{\theta})^2
\end{equation}
with $\tilde{\theta}=\sigma\theta+(1-\sigma)\overline{\theta}$ for some $\sigma\in(0,1)$. It is easy to check that $-1/2<\tilde{\theta}<1/2$ since $|\theta|<1/2$ and $|\overline{\theta}|<1/2$, so that
\[
 m(1-m)\left(\frac{2}{3}\right)^{2-m}\le f''(\tilde{\theta})=\frac{m(1-m)}{(1+\tilde{\theta})^{2-m}}\le  m(1-m)2^{2-m}
\]
so that
\[\begin{split}
{\bf c}f(\overline{\theta}) + {\bf c}\frac{m(1-m)}{\int_\Omega S^{1+m} \dx}\left(\frac{2}{3}\right)^{2-m}\mathcal{E}[\theta]
&\le \overline{\theta}\,'(t) = \frac{{\bf c}\,\int_\Omega  f(\theta)S ^{m+1}\dx}{\int_\Omega S^{1+m} \dx}\\
&\le {\bf c}f(\overline{\theta}) + {\bf c}{\bf c} \frac{m(1-m)}{\int_\Omega S^{1+m} \dx}2^{2-m}\mathcal{E}[\theta] \\
\end{split}
\]
since we recall that
\[
\int_\Omega f'(\overline{\theta})(\theta-\overline{\theta})S^{1+m} \dx=0\qquad\mbox{and}\qquad
\mathcal{E}[\theta]=\frac{1}{2}\int_\Omega\big|\theta(t)-\overline{\theta}(t)\big|^2 S^{1+m} \dx \,.
\]
Finally, we recall that by Proposition \ref{asympt.rate1}, for $t\ge t_0$ we have that
\[
\mathcal{E}[\theta(t)]\le \ee^{-\gamma(t-t_0)}\mathcal{E}[\theta(t_0)]\le \ee^{-\gamma(t-t_0)}\,,
\]
so that
\[\begin{split}
\overline{\theta}\,'(t)
&\le {\bf c}f(\overline{\theta}) + {\bf c} \frac{m(1-m)}{\int_\Omega S^{1+m} \dx}2^{2-m}\mathcal{E}[\theta] \le {\bf c}\left[f(\overline{\theta})+\frac{m(1-m)}{\int_\Omega S^{1+m} \dx}2^{2-m}\ee^{-\gamma(t-t_0)}\right]\\
&:={\bf c}\left[f(\overline{\theta})+k_0\ee^{-\gamma(t-t_0)}\right]
\end{split}
\]
Now define the function $z(t)$ by the relation $f(z)+k_0\ee^{-\gamma(t-t_0)}=0$, i.e. $z(t)=f^{-1}(-k_0\ee^{-\gamma(t-t_0)})$ and $z(t)\to 0$ as $t\to \infty$; it turns out that
\[
z'(t)\ge 0=f(z)+k_0\ee^{-\gamma(t-t_0)}.
\]
Hence $z(t)\ge\overline{\theta}(t)$ for all $t\ge t_0$ whenever $z(t_0)\ge\overline{\theta}(t_0)$. But this fact is in contrast with the fact that $\overline{\theta}(t)\to 0$ as $t\to \infty$: on the line $(t,z(t))$ we have $\overline{\theta}\,'=0$, so we define the regions
\[
Z^+=\left\{(t,y)\;|\; z(t)<y< 0\right\}\qquad\mbox{and}\qquad Z^-=\left\{(t,y)\;|\; y<z(t)\right\}
\]
so that $\overline{\theta}\,'>0$ on $Z^+$ and $\overline{\theta}\,'<0$ on $Z^-$. As a consequence, if $\overline{\theta}(s)\in Z^-$, for some $s\ge t_0$, then $\overline{\theta}(t)\in Z^-$ for all $t\ge s$, since in $Z_-$ we have $\overline{\theta}\,'<0$, therefore $\overline{\theta}(t)$ can not go to zero as $t\to\infty$, which is a contradiction. Finally we have proved that $\overline{\theta}(t)\in Z^+$ for all $t\ge t_0$, which is what we need to conclude that
\[
0\ge\overline{\theta}(t)\ge z(t)=f^{-1}(-k_0\ee^{-\gamma(t-t_0)})
\]
that implies
\[
|\overline{\theta}(t)|\le k_1\ee^{-\gamma(t-t_0)}
\]
since for $|s|<1/2$ we have
\[
\begin{split}
\big|f^{-1}(s)\big|&= \left|f^{-1}(0)+(f^{-1})'(0)s+\frac{1}{2}(f^{-1})''(\tilde{s})s^2\right|\\
&\le \max\left\{|(f^{-1})'(0)|\,,\,|(f^{-1})''(\tilde{s})|\frac{s}{2}\right\}|s|\le k_1 s.
\end{split}
\]
Now we conclude by observing that for all $t\ge t_0$ we have
\[\begin{split}
\frac{1}{2}\int_\Omega\big|\theta(t)\big|^2 S^{1+m} \dx
&\le \int_\Omega\big|\theta(t)-\overline{\theta}(t)\big|^2 S^{1+m} \dx
    +  \int_\Omega \big|\overline{\theta}(t)\big|^2 S^{1+m} \dx
=2 \mathcal{E}[\theta]+ \big|\overline{\theta}(t)\big|^2 \int_\Omega S^{1+m} \dx\\
&\le 2\ee^{-\gamma(t-t_0)}\mathcal{E}[\theta(t_0)]+k_1\ee^{-2\gamma(t-t_0)}\int_\Omega S^{1+m} \dx
\le k_2\ee^{-\gamma(t-t_0)}
\end{split}
\]
which concludes the proof.\qed

\noindent Next we deal with the case $m>1$. Recall that $S$ is the unique stationary state.
\begin{thm}\label{exp.decay.norm.PME}
Under the assumptions of Theorem \ref{asympt.rate.PME}. Then both entropy and the $\LL^2-$norm, decay exponentially with the same rate, more precisely
\begin{equation}\label{rate.entropy.PME}
\int_\Omega\big|\theta(t)\big|^2 S^{1+m}\dx
\le 2\left(\mathcal{E}[\theta(t_1)]\ee^{-(\beta-2)(t-t_1)}+\big|\theta_0\big|^2\int_\Omega S^{1+m}\dx \right)e^{-2(t-t_1)}
\end{equation}
for all $t>t_1\gg 1$, as in Theorem \ref{asympt.rate.PME}.
\end{thm}
\noindent {\sl Proof.~}The first part of the proof is identical to the one of the previous Theorem \ref{exp.decay.norm}, so that we arrive to the differential equation for $\overline{\theta}$
\begin{equation*}
\overline{\theta}\,'(t)=\frac{\rd}{\dt}\frac{\int_\Omega\theta(t,x)S ^{1+m}\dx}{\int_\Omega S^{1+m} \dx}
    = \frac{{\bf c}\,\int_\Omega  f(\theta)S ^{m+1}\dx}{\int_\Omega S^{1+m} \dx}
    \le {\bf c}f(\overline{\theta})\le -{\bf c}(m-1)\overline{\theta}=-\overline{\theta}\,.\\
\end{equation*}
where $f(\theta)=(1+\theta)-(1+\theta)^m$, ${\bf c}=1/(m-1)$, and in the second step we have used the concavity of $f$ together with Jensen inequality, and in the last step we have used that $f$ lies below its tangent at the origin
\[
f(\overline{\theta})\le -(m-1)\overline{\theta}\,.
\]
An integration over $[s,t]\subset [t_1,\infty)$ gives
\[
\overline{\theta}(t)\le \overline{\theta}(s)\ee^{-(t-s)}\qquad\mbox{and}\qquad
\big|\overline{\theta}(t)\big|\le \big|\overline{\theta}(s)\big|\ee^{-(t-s)}
\]
where $t_1$ is as in Theorem \ref{asympt.rate.PME}. Moreover we have that
\begin{equation*}
\frac{\mathcal{E}[\theta(t)]}{\int_\Omega S^{1+m}\dx}=\frac{1}{2}\int_\Omega\big|\theta(t)-\overline{\theta}(t)\big|^2
    \frac{S^{1+m}\dx}{\int_\Omega S^{1+m}\dx}
    =\frac{\int_\Omega\big|\theta(t)\big|^2 S^{1+m}\dx}{2\int_\Omega S^{1+m}\dx}
    -\frac{1}{2}\big|\overline{\theta}(t)\big|^2
\end{equation*}
Combining the above result with the entropy decay of Theorem \ref{asympt.rate.PME}
\begin{equation*}
\mathcal{E}[\theta(t)]\le \mathcal{E}[\theta(t_1)]\ee^{-\beta(t-t_1)}\qquad\mbox{for all }t\ge t_1.
\end{equation*}
we obtain
\[\begin{split}
\int_\Omega\big|\theta(t)\big|^2 S^{1+m}\dx
&=2\mathcal{E}[\theta(t)]+2\big|\overline{\theta}(t)\big|^2\,\int_\Omega S^{1+m}\dx\\
&\le 2\left(\mathcal{E}[\theta(t_1)]\ee^{-(\beta-2)(t-t_1)}+\big|\overline{\theta}_0\big|^2\int_\Omega S^{1+m}\dx \right)e^{-2(t-t_1)}\,\mbox{.\qed}
\end{split}\]

\section{Stationary solutions and their limit as $p\to 1$.}\label{ell-sect}

Let $1\le p<p_s$  and let $U_p$ be a solution to the elliptic problem
\begin{equation}\label{ell.prob.p}
\left\{
\begin{array}{lll}
-\Delta U=\lambda_p\,U^p&\mbox{in }\Omega\\
U>0&\mbox{in }\Omega\\
U=0&\mbox{on }\partial\Omega\\
\end{array}
\right.
\end{equation}
where $\lambda_p>0$ if $1<p<p_s$ and $\lambda_p=\lambda_1$ for $p=1$. We are interested in the relation between solutions of the elliptic equation for different values of $p\in[1, p_s)$, in particular we would like to see whether the limit $V:=\lim_{p\to 1}U_p$ exists and under which conditions it is the ground state of the Dirichlet Laplacian $\Phi_1$ on $\Omega$. The existence of a limit depends on a normalization that we will discuss below.

It is well understood by subcritical semilinear theory that positive weak solutions of the above elliptic problem are indeed classical solutions up to the boundary. Weak solutions can be defined as follows: a function $U_p\in W_0^{1,2}(\Omega) $ is a weak solution to the elliptic problem \eqref{ell.prob.p} if and only if
\begin{equation}\label{weak.ell.p}
\int_\Omega\left[\nabla U_p\cdot\nabla\varphi-\lambda_p U_p^p\varphi\right]\dx=0
\end{equation}
for all $\varphi\in W_0^{1,2}(\Omega)$. Notice that when $p=1$ there is a positive solution, unique up to a multiplicative constant, while when $p>1$ uniqueness is not always true, it depends on the geometry of the domain. The difficulty in understanding the limit of $U_p$ as $p\to 1^+$, relies indeed in the lack of uniqueness and on a scaling property typical of the nonlinear problem. In the case of uniqueness, for example in the case when $\Omega$ is a ball, solutions are variational, in the sense that they are minima of a the functional $\|\nabla U\|_2^2$ under the restriction $\|U\|_{p+1}=1$, but when the uniqueness is not guaranteed, solutions are just critical points of such functional.

One can also easily see that the constant $\lambda_p>0$ in the nonlinear problem can be manipulated by rescaling, because if $U_{p,(1)}(x)$ is a
solution with parameter $\lambda_{p,(1)}$, then $U_{p,(2)}(x)=\mu^{1/(p-1)}\,U_{p,(2)}(x)$ is a solution with parameter $\lambda_{p,(2)}=\mu \lambda_{p,(1)}$. In any normed space $\|U_{p,(2)}\|=\mu^{1/(p-1)}\|U_{p,(1)}\|.$  This means that scaling allows to fix the norm of a solution: changing the norm by a factor $\mu^{1/(p-1)}$ by  scaling is equivalent to changing $\lambda_p$ in the equation by a factor $\mu^{-1}$.

\noindent\textbf{Assumption throughout this section.} Let us fix $\lambda_p$ as the factor for which $\|U_p\|_{p+1}=1$, so that, using $U_p$ as test function, we obtain the following identity
\begin{equation}\label{energy.id}
\|\nabla U_p\|_2^2=\lambda_p\|U_p\|_{p+1}^{p+1}=\lambda_p\,,
\end{equation}
so that it is  equivalent to prove that $\lambda_p\to \lambda_1$ or to prove that $\|\nabla U_p\|_2\to \|\nabla\Phi_1\|_2$, when $p\to 1$. Recall that $\Phi_1$ has unit $\LL^2$-norm.

We state now the main result of this section.
\begin{thm}\label{Elliptic.THM}
Let $U_p$ be a family of solutions of Problem {\rm \ref{ell.prob.p}} with $p\in[1,p_s)$, $\|U_p\|_{p+1}=1$ and let $\lambda_p>0$ be chosen according to \eqref{energy.id}. Then as $p\to 1$, $\lambda_p\to \lambda_1$, $U_p\to \Phi_1$ in $\LL^\infty(\Omega)$\,,\, $\nabla U_p\to\nabla \Phi_1$ in $\big(\LL^2(\Omega)\big)^d$. Besides, there exist two explicit constants $0<c_0<c_1$ such that
\begin{equation}
c_0^{p-1}\lambda_1\le\lambda_p\le c_1^{p-1}\lambda_1\,.
\end{equation}
Moreover, there exists constants $0< \widetilde{k}_0(p)\le \widetilde{k}_1(p)$ such that $\widetilde{k}_i(p)\to 1$ as $p\to 1^+$, such that
\begin{equation}
\widetilde{k}_0(p) \le \frac{U_p(x)}{\Phi_1(x)}\le \widetilde{k}_1(p),\qquad\mbox{for all }x\in \overline{\Omega}\,.
\end{equation}
\end{thm}
\subsection{Proof of Theorem \ref{Elliptic.THM}}
The proof of the above Theorem will be divided into several steps. Our first result in this connection is the following.
\begin{lem}\label{conv.lambda.p}
Let $U_p$ be a solution of Problem {\rm \ref{ell.prob.p}} with $p\in[1,p_s)$, $\|U_p\|_{p+1}=1$ and let $\lambda_p>0$ be chosen according to \eqref{energy.id}. If there is a constant $A>0$ such that
\[
0<\lambda_p \le A<\infty\,,
\]
then $U_p\to \Phi_1$ in $\LL^q(\Omega)$ for any $0<q<2^*$, and $\lambda_p\to \lambda_1$.
\end{lem}
\noindent {\sl Proof.~}Since $U_p$ is a solution to the elliptic Problem \ref{ell.prob.p} with $\|U_p\|_{p+1}=1$ we have that the hypotheses together with the energy identity \eqref{energy.id} give
\[
0<\lambda_p=\|\nabla U_p\|_2^2 \le A <\infty
\]
which proves that $\lambda_p=\|\nabla U_p\|_2^2$ is uniformly bounded for all $p\in[1,p_s)$. Hence we can guarantee that there exists a subsequence $\nabla U_{p_n}$ that converges weakly in $\LL^2(\Omega)$ to a function $W$. Moreover, by Kondrachov's compactness theorem  there is a (maybe different) subsequence $U_{p_k}$ that converges to $V$ strongly in any $\LL^q$ with $1\le q<2^*$. Strong convergence implies that $\|V\|_2=1$, hence $V$ can not be identically zero. Moreover, it is well known that in this case $W=\nabla V$.
Next, we show that $U_{p_n}^{p_n}\to V$ in $L^1(\Omega)$:
\[\begin{split}
\int_\Omega \left|U_{p_n}^{p_n}-V\right|\dx
    &\le \int_\Omega \left|U_{p_n}^{p_n}-U_{p_n}\right|\dx
        +\int_\Omega \left|U_{p_n}-V\right|\dx
    = \int_\Omega \left|U_{p_n}^{p_n-1}-1\right| U_{p_n}\dx
        +\int_\Omega \left|U_{p_n}-V\right|\dx
\end{split}
\]
where the second integral converges to zero since $U_{p_n}\to V$ in $\LL^1(\Omega)$, while for the first we have to use a numerical inequality:
\[
\left|a^z-1\right|\le \left(\frac{a^b}{b}+\frac{1}{a}\right)\,z\qquad\mbox{for all}\qquad a>0\,,~~0\le z\le b
\]
which we will prove at the end of the proof. Using the above numerical inequality in the first step for $a=U_{p_n}>0$ and $0\le z=p_n-1\le b=2^*-2=4/(d-2)$, we obtain
\[\begin{split}
\int_\Omega \left|U_{p_n}^{p_n-1}-1\right| U_{p_n}\dx
    &\le (p_n-1)\int_\Omega \left(\frac{U_{p_n}^b}{b}+\frac{1}{U_{p_n}}\right) U_{p_n}\dx
     \le \frac{p_n-1}{b}\int_\Omega \left(U_{p_n}^{b+1}+b\right) \dx\\
    &\le \frac{p_n-1}{b}\left[b|\Omega|+|\Omega|^{1-\frac{b+1}{2^*}}
        \left(\int_\Omega U_{p_n}^{2^*} \dx\right)^{\frac{b+1}{2^*}}\right]\\
    &\le \frac{p_n-1}{b}\left[2^*|\Omega|+|\Omega|^{1-\frac{b+1}{2^*}}
        \left(\int_\Omega |\nabla U_{p_n}|^2 \dx\right)^{\frac{b+1}{2}}\right]\\
    &= \frac{p_n-1}{b}\left[2^*|\Omega|+|\Omega|^{1-\frac{b+1}{2^*}}
        \left(\mathcal{S}_2\lambda_p\right)^{\frac{b+1}{2}}\right]
     \le \frac{p_n-1}{b}\left[2^*|\Omega|+|\Omega|^{1-\frac{b+1}{2^*}}
        \left(\mathcal{S}_2\,A\right)^{\frac{b+1}{2}}\right]
\end{split}
\]
since $b+1\le 2^*-1=p_s=(d+2)/(d-2)$, we can use Sobolev and H\"older inequalities, and the fact that the $\LL^2$-norm of the gradient, or equivalently $\lambda_p$, is uniformly bounded by $A$.

\noindent Now, we identify the limits. For all $1\le p<p_s$:
\begin{equation}\label{eq.lambda.1.p}
\lambda_p\int_\Omega U_p^p\,\Phi_1 \dx=-\int_\Omega \Phi_1 \Delta U_p \dx =\int_\Omega \nabla U_p\cdot\nabla\Phi_1 \dx =-\int_\Omega U_p\Delta\Phi_1 \dx =\lambda_1\int_\Omega U_p\,\Phi_1 \dx
\end{equation}
such equalities hold by the weak form of the equation satisfied by $U_p$ and since both $U_p$ and $\Phi_1$ are in $W_0^{1,2}(\Omega)$. Take any subsequence of $p_n$ such that $\lambda_{p_n}$ converges. Let $\Lambda=\lim_{k\to\infty}\lambda_{p_{n_k}}$. Taking limits as $k\to \infty$ so that $p_{n_k}\to 1$ in the above expression to get
\[
\Lambda\int_\Omega V\,\Phi_1 \dx=\lambda_1\int_\Omega V\,\Phi_1 \dx
\]
since we know that both $U_{p_{n_k}}^{p_{n_k}-1}$ and $U_{p_{n_k}}$ converge to $V$. We conclude that $\lim_{k\to\infty}\lambda_{p_{n_k}}=\lambda_1$. Since this holds for all subsequences, this means that $\lim_{n\to\infty}\lambda_{p_n}=\lambda_1$. We are now ready to identify $V=\lim_{n\to \infty}U_{p_n}$, indeed we just notice that
\[
\int_\Omega\nabla V\cdot\nabla\varphi\dx-\lambda_1 \int_\Omega V\varphi\dx = \lim_{n\to \infty} \left[\int_\Omega\nabla U_{p_n}\cdot\nabla\varphi\dx -\lambda_{p_n}\int_\Omega U_{p_n}^{p_n}\varphi\dx\right]=0
\]
where all the quantities have been shown to converge in such a way. The latter equality identifies $V$ as the unique ground state $\Phi_1$ such that $-\Delta \Phi_1=\lambda_1\Phi_1$ and $\|V\|_2=\|\Phi_1\|_2=1$.

Now we prove that for any subsequence $p_n\to 1$, we have $\lambda_{p_n}\to \lambda_1$, $U_{p_n}\to \Phi_1$ and consequently $U_{p_n}^{p_n}\to \Phi_1$. Suppose that there exists a sequence $p_n$ such that $\lim_n\lambda_{p_n}=\Lambda\neq \lambda_1$. We can repeat the first steps to conclude that there is a subsequence $U_{p_{n_k}}$ that converges strongly in $L^1$ to some $V$ and such that $\nabla U_{p_{n_k}}$ converges weakly in $\LL^2$ to $W=\nabla V$. Moreover, we also have that $U_{p_n}^{p_n}\to V$. Using formula \eqref{eq.lambda.1.p}, we have that
\[
\lambda_{p_{n_k}} \int_\Omega U_{p_{n_k}}^{p_{n_k}}\,\Phi_1 \dx = \lambda_1 \int_\Omega U_{p_{n_k}}\,\Phi_1 \dx
\]
and taking the limit as $k\to\infty$ we get a contradiction, since $\lim_k\lambda_{p_{n_k}}=\lambda_1$. The proof is concluded once we prove the numerical inequality
\[
\left|a^z-1\right|\le \max\left\{\frac{a^b-1}{b},\frac{1}{a}\right\}\le\left(\frac{a^b}{b}+\frac{1}{a}\right)\,z\qquad\mbox{for all}\qquad a>0\,,~~0\le z\le b.
\]
we prove it first for $a>1$: since the function $f(z)=a^z-1$ is convex, it lies below the secant $(a^b-1)z/b$ for all $0\le z\le b$, hence the inequality $|a^z-1|=a^z-1\le (a^b-1)z/b$ when $a>1$ and $0\le z\le b$. When $0<a<1$, we see that $|a^z-1|=1-a^z$. We see that $f(z)=1-a^z=z/a$ when $z=0$, and that $f'(z)=-\log(a)a^z\le \log(1/a)\le 1/a$, hence the desired inequality is valid also when $0<a<1$.\qed

\medskip

\noindent In view of this result, we need  to prove the upper bound $\lambda_p\le A$.
The proof is based on an idea of  Brezis and Turner \cite{BT}. This method relies on an Hardy-type inequality, which holds for a large class of domains, but only in the range $1\le p\le (d+1)/(d-1)$. If one wants to deal with the full range of exponents $1\le p<p_s$, one has to proceed as Gidas-Ni-Nirenberg \cite{GNN} when the domain is convex, or as DeFigueredo-Lions-Nussbaum \cite{dFLN} which extend the ideas of \cite{GNN} to more general domains.

\begin{prop}The following Hardy-type inequality holds true whenever $\Omega$ has a finite inradius and satisfies a uniform exterior ball condition
\begin{equation}\label{Hardy.BT}
\left\|\frac{f}{\Phi_1^r} \right\|_q\le H_{r,d}\|\nabla f\|_2\qquad\mbox{if }f\in W^{1,2}_0(\Omega),~0<q\le\frac{2d}{d-2+2r}\,,~\mbox{and } 0\le r\le 1\,.
\end{equation}
where $\Phi_1$ is the unique positive ground state of the Dirichlet Laplacian on $\Omega$, and $H_{r,d}$ is a suitable positive constant that depends only on $r,d$ and $|\Omega|$ and is given at the end of the proof.
\end{prop}
\noindent {\sl Proof.~} The proof is obtained by combining the standard Hardy inequality
\[
\left\|\frac{f}{\dist(\cdot, \partial\Omega)} \right\|_2\le H_0\|\nabla f\|_2\qquad\mbox{for all }f\in W^{1,2}_0(\Omega)\,,
\]
which holds whenever $\Omega$ has a finite inradius and satisfies a uniform exterior ball condition, see for instance Section 1.5 of \cite{Davies}. We combine such Hardy inequality with the standard Sobolev imbedding $\|f\|_{2^*}\le \mathcal{S}_2^2\|\nabla f\|_2^2$ as follows:
\[\begin{split}
\int_\Omega \frac{f^q}{\dist(x,\partial\Omega)^rq}\dx
&=\int_\Omega\frac{f^{rq}}{\dist(x,\partial\Omega)^{rq}} u^{(1-r)q}\dx\\
&\le \left[\int_\Omega\left(\frac{f^{rq}}{\dist(x,\partial\Omega)^{rq}}\right)^{\gamma}\dx\right]^{\frac{1}{\gamma}}
\left[\int_\Omega\left(f^{(1-r)q}\right)^{\frac{\gamma}{\gamma-1}}\dx\right]^{\frac{\gamma-1}{\gamma}}\\
&=_{(a)}\left[\int_\Omega\frac{u^2}{\dist(x,\partial\Omega)^2}\dx\right]^{\frac{rq}{2}}
\left[\int_\Omega f^{\frac{2(1-r)q}{2-rq}}\dx\right]^{\frac{2-rq}{2}}\\
&\le_{(b)}\left[H_0\|\nabla f\|_2\right]^{rq}|\Omega|^{1-q\frac{d-2+2r}{2d}}\mathcal{S}_2^{(1-r)q}\|\nabla f\|_2^{(1-r)q}
\end{split}
\]
where in $(a)$ we have set  $\gamma=2/rq>1$ and $\gamma/(\gamma-1)=2/(2-rq)$, for any $0\le r\le 1$, while in $(b)$ we have used the above Hardy inequality and we have estimated the second integral with the H\"older inequality
\[\begin{split}
\left[\int_\Omega f^{\frac{2(1-r)q}{2-rq}}\dx\right]^{\frac{2-rq}{2}}
&=\|f\|_\frac{2(1-r)q}{2-rq}^{(1-r)q}
\le |\Omega|^{(1-r)q\left(\frac{2-rq}{2(1-r)q}-\frac{1}{2^*}\right)}\|f\|_{2^*}^{(1-r)q}\\
&\le |\Omega|^{1-q\frac{d-2+2r}{2d}}\mathcal{S}_2^{(1-r)q}\|\nabla f\|_2^{(1-r)q}
\end{split}\]
and in order to use Sobolev inequality, we need
\[
\frac{(1-r)q\gamma}{\gamma-1}= \frac{2(1-r)q}{2-rq}\le 2^*\qquad\mbox{that is}\qquad q\le \frac{2d}{d-2+2r}\,.
\]
We have obtained so far
\begin{equation}\label{Hardy.dist}
\left[\int_\Omega \frac{u^q}{\dist(x,\partial\Omega)^{rq}}\dx\right]^{\frac{1}{q}}\le H_0^{r}|\Omega|^{\frac{1}{q}-\frac{d-2+2r}{2d}}\mathcal{S}_2^{1-r}\|\nabla f\|_2
\end{equation}
We conclude the proof by noticing that, since there exists two positive constants $c_0,c_1$ depending only on the dimension $d\ge 3$, such that
\[
c_0\,\dist(\cdot, \partial\Omega)\le \Phi_1(\cdot) \le c_1\,\dist(\cdot, \partial\Omega)\qquad\mbox{in }\Omega\,.
\]
We combine the above lower bound together with \eqref{Hardy.dist} to get the desired inequality \eqref{Hardy.BT}
\[
\left[\int_\Omega \frac{u^q}{\Phi_1^{rq}}\dx\right]^{\frac{1}{q}}\le \frac{1}{c_0^r}\left[\int_\Omega \frac{u^q}{\dist(x,\partial\Omega)^{rq}}\dx\right]^{\frac{1}{q}}\le \frac{H_0^{r}}{c_0^r}|\Omega|^{\frac{1}{q}-\frac{d-2+2r}{2d}}\mathcal{S}_2^{1-r}\|\nabla f\|_2:=H_{r,d}\|\nabla f\|_2\,\mbox{.\qed}
\]
We are now ready to prove the upper bounds for $\lambda_p$.
\begin{prop}\label{upper.lambda.p.1}
Let $1\le p\le(d+1)/(d-1)$ and $d\ge 3$ and $\lambda_p$ be such that $\|U_p\|_{p+1}=1$, as in \eqref{energy.id}. Then the following upper bound holds true
\begin{equation}\label{upper.lambda.p.BT}
\lambda_p^{\frac{[d+1-p(d-1)](p+1)}{2(p-1)}}\le\lambda_1^{\frac{2p}{p-1}} \left(
    \int_\Omega\Phi_1\dx\right)^2 \, H_{r,d}^{(d-1)p+d+1}\,.
\end{equation}
\end{prop}
\noindent {\sl Proof.~}Testing the equation \eqref{ell.prob.p} with $\Phi_1$ yields, as in \eqref{eq.lambda.1.p}
\[
\lambda_p\int_\Omega U_p^p\Phi_1\dx=\lambda_1\int_\Omega U_p\Phi_1\dx\qquad\mbox{or}\qquad \frac{\lambda_p}{\lambda_1}\int_\Omega U_p^p\Phi_1\dx=\int_\Omega U_p\Phi_1\dx
\]
that gives
\[
\frac{\lambda_p}{\lambda_1}\int_\Omega U_p^p\Phi_1\dx=\int_\Omega U_p\Phi_1\dx
\le \left(\int_\Omega U_p^p\Phi_1\dx\right)^{\frac{1}{p}}\left(\int_\Omega\Phi_1\dx\right)^{1-\frac{1}{p}} \,.
\]
where we have used H\"older inequality in the last step. We have obtained that
\begin{equation}\label{U.p.p}
\int_\Omega U_p^p\Phi_1\dx \le\left(\frac{\lambda_1}{\lambda_p}\right)^{\frac{p}{p-1}}\int_\Omega\Phi_1\dx\,.
\end{equation}

Next we calculate
\[\begin{split}
1&=\int_\Omega U_p^{p+1}\dx
 =\int_\Omega\left(U_p^p\Phi_1\right)^{\alpha} \left(\frac{U_p^q}{\Phi_1^{qr}}\right)^{1-\alpha}\dx
\le_{(a)} \left[\int_\Omega U_p^p\Phi_1\dx\right]^{\alpha} \left[\int_\Omega \frac{U_p^q}{\Phi_1^{qr}}\right]^{1-\alpha}\\
&\le_{(b)} \left[ \left(\frac{\lambda_1}{\lambda_p}\right)^{\frac{p}{p-1}}
    \int_\Omega\Phi_1\dx\right]^{\alpha} \,\left\|\frac{U_p}{\Phi_1^r} \right\|_q^{q(1-\alpha)}
\le_{(c)} \left[ \left(\frac{\lambda_1}{\lambda_p}\right)^{\frac{p}{p-1}}
    \int_\Omega\Phi_1\dx\right]^{\alpha} \, \left(H_{r,d}\|\nabla U_p\|_2\right)^{q(1-\alpha)}\\
&=_{(d)} \left[ \left(\frac{\lambda_1}{\lambda_p}\right)^{\frac{p}{p-1}}
    \int_\Omega\Phi_1\dx\right]^{\alpha} \, \left(H_{r,d}\lambda_p^{\frac{1}{2}}\right)^{(1-\alpha)p+1}
\end{split}
\]
where in $(a)$ we have used H\"older inequality with conjugate exponents $1/\alpha$ and $\alpha/(1-\alpha)$, for some $\alpha\in(0,1)$ to be fixed later. We have also put
\[
p+1=\alpha p+(1-\alpha)q\qquad\mbox{and}\qquad \alpha=r(1-\alpha)q
\]
which is equivalent to fix the values of $q$ and $r$ (as functions of $1\le p<p_s$ and $0<\alpha<1$) as follows
\[
q=p+\frac{1}{1-\alpha}\qquad\mbox{and}\qquad r=\frac{1}{(1-\alpha)q}=\frac{\alpha}{1+(1-\alpha)p}
\]
we notice that $0\le r\le 1$ since $0<\alpha<1$ and $p\ge 1$. In $(b)$ we have used the upper bounds \eqref{U.p.p}\,, while in $(c)$ we have used the Hardy inequality \eqref{Hardy.BT}, for which we have to check that $q\le 2d/(d-2+2r)$, that is equivalent to
\[
p+\frac{1}{1-\alpha}\le \frac{2d}{d-2+\frac{2\alpha}{1+(1-\alpha)p}}\qquad\mbox{that is}\qquad \alpha\le\frac{d+2-(d-2)p}{2(d+1)-(d-2)p}<1
\]
since $p<p_s=(d+2)/(d-2)$ and hence $0<d+2-(d-2)p<2(d+1)-(d-2)p$. Finally, in the last step $(d)$ we have used the identity $\|\nabla U_p\|_2^2=\lambda_p\|U_p\|_{p+1}^{p+1}=\lambda_p$ and the fact that $q(1-\alpha)=(1-\alpha)p+1$.

We have obtained the desired upper bound \eqref{upper.lambda.p.BT}
\[
\lambda_p^{\frac{p\alpha}{p-1}-\frac{(1-\alpha)p+1}{2}}\le\lambda_1^{\frac{\alpha p}{p-1}} \left(
    \int_\Omega\Phi_1\dx\right)^{\alpha} \, H_{r,d}^{(1-\alpha)p+1}\,,
\]
for any $\alpha\in(0,1)$ such that
\[
\frac{p\alpha}{p-1}-\frac{(1-\alpha)p+1}{2}\ge 0\qquad\mbox{that is}\qquad \alpha\ge\frac{p-1}{p}
\]
that is
\[
\frac{p-1}{p}\le \alpha\le \frac{d+2-(d-2)p}{2(d+1)-(d-2)p}
\]
which is non empty only when $p\le (d+1)/(d-1)$. Letting now $\alpha=2/(d+1)$ gives the desired upper bound \eqref{upper.lambda.p.BT}, since
\[
\frac{p\alpha}{p-1}-\frac{(1-\alpha)p+1}{2}=\frac{-(d-1)p^2+2p+d+1}{2(d+1)(p-1)}
=\frac{\big[d+1-p(d-1) \big](p+1)}{2(d+1)(p-1)}\mbox{\,.\qed}
\]
At this point we are able to prove the first part of Theorem \ref{Elliptic.THM}.
\begin{prop}\label{prop.47}
Let $U_p$ be a family of solution of Problem {\rm \ref{ell.prob.p}} with $p\in[1,p_s)$, $\|U_p\|_{p+1}=1$ and let $\lambda_p>0$ be chosen according to \eqref{energy.id}. Then as $p\to 1$, $\lambda_p\to \lambda_1$, $U_p\to \Phi_1$ in $\LL^\infty(\Omega)$\,,\, $\nabla U_p\to\nabla \Phi_1$ in $\big(\LL^2(\Omega)\big)^d$. Moreover, there exist two explicit constants $0<c_0<c_1$ such that
\begin{equation}
c_0^{p-1}\lambda_1\le\lambda_p\le c_1^{p-1}\lambda_1\,.
\end{equation}
\end{prop}
\noindent {\sl Proof.~}The only thing that remains to prove is the convergence in $\LL^\infty$ and the convergence of the gradients. We know that $U_p$ and $\Phi_1$ are H\"older continuous on the whole $\overline{\Omega}$, and the $C^\alpha$-norm of both functions is uniformly bounded, say $\|U_p\|_{C^\alpha(\Omega)}+\|\Phi_1\|_{C^\alpha(\Omega)}\le K$. Moreover combining the upper estimates for $\lambda_p$ of Proposition \ref{upper.lambda.p.1} (or of Proposition \ref{lambda.prop} when there is uniqueness of the stationary state) with Lemma \ref{conv.lambda.p} gives the convergence of $U_p\to \Phi_1$ in $\LL^q(\Omega)$ for any $1\le q<2^*$. By the interpolation Lemma \ref{GN}, we get that $U_p\to \Phi_1$ in $\LL^\infty(\Omega)$, since when $p\to 1$ we have
\[
\|U_p-\Phi_1\|_\infty\le C\,\|U_p-\Phi_1\|_{C^\alpha(\Omega)}^{1-\vartheta}\|U_p-\Phi_1\|_2^{\vartheta}
\le C\,K^{1-\theta}\|U_p-\Phi_1\|_2^{\theta}\to 0\,.
\]
As for the gradients, we use the Hardy inequality \eqref{Hardy.BT} with $r=1$ and $q=2$ applied to $U_p-\Phi_1\in W_0^{1,2}(\Omega)$, to get
\[
\int_\Omega \frac{|U_p-\Phi_1|^2}{\Phi_1^2}\dx \le H_{1,2}\int_\Omega \big|\nabla (U_p-\Phi_1)\big|^2\dx
\]
Next, we analyze the right-hand side:
\[
\begin{split}
&\int_\Omega \big|\nabla (U_p-\Phi_1)\big|^2\dx
= \int_\Omega \big|\nabla U_p\big|^2\dx + \int_\Omega \big|\nabla \Phi_1\big|^2\dx
    -2\int_\Omega \nabla U_p\cdot\nabla \Phi_1\dx\\
(a)&= \lambda_p\int_\Omega U_p^{p+1}\dx+\lambda_1\int_\Omega\Phi_1^2\dx
- \lambda_p\int_\Omega U_p^p\,\Phi_1 \dx -\lambda_1\int_\Omega U_p\,\Phi_1\dx\\
&\le\lambda_p\int_\Omega U_p^p\big|U_p-\Phi_1\big|\dx + \lambda_1\int_\Omega \big|U_p-\Phi_1\big|\,\Phi_1\dx\\
(b)&\le\lambda_p\left[\int_\Omega U_p^{p+1}\dx\right]^{\frac{p}{p+1}} \left[\int_\Omega\big|U_p-\Phi_1\big|^{p+1}\dx\right]^{\frac{1}{p+1}}
+ \lambda_1 \left[\int_\Omega \big|U_p-\Phi_1\big|^2\dx\right]^{\frac{1}{2}}\,
\left[\int_\Omega\Phi_1^2\dx\right]^{\frac{1}{2}},\\
\end{split}
\]
where in $(a)$ we have used formula \eqref{eq.lambda.1.p} in the form
\[
2\int_\Omega \nabla U_p\cdot\nabla \Phi_1\dx
=\lambda_p\int_\Omega U_p^p\,\Phi_1 \dx+\lambda_1\int_\Omega U_p\,\Phi_1 \dx\,.
\]
As for $(b)$ we have used H\"older inequality with conjugate exponents $(p+1)/p$ and $p+1$ for the first term, and Cauchy-Schwartz inequality for the second term. Putting all the pieces together we have obtained
\[
\left\|\frac{U_p-\Phi_1}{\Phi_1}\right\|_2^2
\le\lambda_p \|U_p\|_{p+1}^p \|U_p-\Phi_1\|_{p+1}
+ \lambda_1 \|U_p-\Phi_1\|_2\|\Phi_1\|_2\;\xrightarrow[p\to\,1^+]{} 0
\]
since we already know by Lemma \ref{conv.lambda.p} that $U_p\to\Phi_1$ in $\LL^q(\Omega)$ for any $1\le q <2^*$, and we also know that $p+1<2^*$ since $p<p_s$. The sole requirement of Lemma \ref{conv.lambda.p} is that $\lambda_p\le A$, and the uniform upper bound for $\lambda_p$ is guaranteed by Proposition \ref{upper.lambda.p.1} when $1\le p\le(d+1)/(d-1)$ or by Proposition \ref{lambda.prop} when $U_p$ is a variational solution and $1\le p<p_s$.\qed

The last step in proving Theorem \ref{Elliptic.THM} consists in comparing solutions corresponding to different $p$ and $\lambda_p$, more precisely to show that there exist constants $0< \widetilde{k}_0(p)\le \widetilde{k}_1(p)$ such that $\widetilde{k}_i(p)\to 1$ as $p\to 1^+$. $\Phi_1$ is the corresponding ground state, towards to $U_p$ converges as $p\to 1$.
\begin{prop}\label{prop.k0.k1}
Under the running assumptions on $U_p$ and $\Phi_1$, there exist constants $0< \widetilde{k}_0(p)\le \widetilde{k}_1(p)$ such that $\widetilde{k}_i(p)\to 1$ as $p\to 1^+$, such that
\begin{equation}
\widetilde{k}_0(p) \le \frac{U_p(x)}{\Phi_1(x)}\le \widetilde{k}_1(p),\qquad\mbox{for all }x\in \overline{\Omega}\,.
\end{equation}
\end{prop}
\noindent {\sl Proof.~}The proof is divided in several steps.

\noindent$\bullet~$\textsc{Step 1. }\textit{Convergence of the quotient in an inner region. }Proposition \ref{prop.47} implies that $U_p/\Phi_1\to 1$ in any inner region in which $\Phi_1\ge \sigma>0$. In the sequel we will construct a special region as follows.
By Lemma \ref{Lemma.Distance} we know that there exists a $\delta>0$ such that
\[
V_{\delta}=\left\{x\in\Omega\,:\, d(x)<\delta \right\}
\]
such that $d(x)=\dist(x,\partial\Omega)\in C^2(\Omega_{\delta})$ and $d(x)$ is Lipschitz with constant
$1$, i.e. $|d(x)-d(y)|\le\, |x-y|$, and $0<c\le\,|\nabla d(x)|\le 1$, $-K\le\Delta d(x)\le K$ in $\Omega_\delta$. In the complement, $\Omega_\delta=\Omega-V_\delta$, we know that
\[
\sigma\le U_p(x), \Phi_1(x) \le M\qquad\mbox{for all }x\in\Omega_\delta,
\]
so that, given $\varepsilon>0$ for $p$ sufficiently close to $1$ we have
\[
(1-\varepsilon)\Phi_1(x)\le U_p(x)\le(1+\varepsilon)\Phi_1(x)\qquad\mbox{for all $x\in\Omega_\delta$}\,.
\]
It remains to prove that the above inequality extends to the thin region $V_\delta=\Omega-\Omega_\delta$, and this will be done in the next steps.

\noindent$\bullet~$\textsc{Step 2. }\textit{Upper comparison near the boundary.} The upper estimate for $\lambda_p$ of Theorem \ref{Elliptic.THM} reads $\lambda_p\le c_1^{p-1}\lambda_1$. Since we are working in the thin annular domain $V_\delta=\Omega\setminus\Omega_\delta$, and we know that $\Phi_1=0$ on $\partial\Omega$, we can assume (eventually by taking a smaller $\delta>0$) that $(c_1(1+\varepsilon)\Phi_1)^{p-1}\le 1$ in $V_\delta$. As a consequence, we have that
\[
\lambda_p\big((1+\varepsilon)\Phi_1\big)^p \le\lambda_1 \big(c_1(1+\varepsilon)\Phi_1\big)^{p-1} (1+\varepsilon)\Phi_1\le \lambda_1(1+\varepsilon)\Phi_1\,.
\]
This allows to compare $U_p$ and $\Phi_2=(1+\varepsilon)\Phi_1$ on the thin set $V_\delta$. We have $0\le U_p\,, (1+\varepsilon)\Phi_1\le M$ in $\overline{V_\delta}$. The respective equations are
\[\left\{
\begin{array}{lll}
-\Delta U_p=\lambda_p\,U_p^p\,, &\mbox{in }V_\delta\,,\\
-\Delta \Phi_2=-\Delta (1+\varepsilon)\Phi_1 =\lambda_1 (1+\varepsilon)\Phi_1
    \ge \lambda_p\big((1+\varepsilon)\Phi_1\big)^p=\lambda_p\Phi_2^p \, &\mbox{in }V_\delta\,,\\
\end{array}\right.
\]
and the boundary data
\[\left\{
\begin{array}{lll}
\Phi_2 = U_p =0\,, & \mbox{on }\partial\Omega\,,\\
\Phi_2 \ge U_p \ge 0\,, & \mbox{on }\partial\Omega_\delta\,.\\
\end{array}\right.
\]
We want to apply Theorem \ref{thm.upper.comp} to obtain the comparison  $U_p\le (1+\varepsilon)\Phi_1=\Phi_2$ in $V_\delta$. We need a the following smallness condition on $V_\delta$:
\[
|V_\delta|<\frac{\omega_d}{\left(2p\,\lambda_p\,M^{p-1}\right)^d}\,.
\]
The above condition can be fulfilled just by choosing $\delta$ sufficiently small, and we can always do that independently of $\varepsilon$ small.

\noindent$\bullet~$\textsc{Step 3. }\textit{Lower comparison near the boundary. }This part consists of two comparison arguments. First we observe that we can compare $U_p$ with a suitable harmonic function on $V_\delta$, namely
\[\left\{
\begin{array}{lll}
-\Delta U_p=\lambda_p\,U_p^p\ge 0\,, &\mbox{in }V_\delta\,,\\
-\Delta U= 0&\mbox{in }V_\delta\,,\\
\end{array}\right.
\]
and the boundary data
\[\left\{
\begin{array}{lll}
U = U_p =0\,, & \mbox{on }\partial\Omega\,,\\
U=(1-\varepsilon)\Phi_1 \le U_p \,, & \mbox{on }\partial\Omega_\delta\,.\\
\end{array}\right.
\]
we apply standard comparison to get $U\le U_p$ in $V_\delta$.

Next we want to prove that $U\ge (1-2\varepsilon)\Phi_1$ on $V_\delta$ if $\delta$ is small enough.

We define the function $w=(1-\varepsilon)\Phi_1-U$ that satisfies
\[\left\{
\begin{array}{lll}
-\Delta w=\lambda_1(1-\varepsilon)\Phi_1\le c_0\delta \,, &\mbox{in }V_\delta\,,\\
w= 0\,, &\mbox{on }\partial V_\delta\,,\\
\end{array}\right.
\]
and we compare $\hat{w}=w/(c_0\delta)$ with $W$, which solves the following problem on the whole $\Omega=\overline{\Omega_\delta}\cup V_\delta$:
\[\left\{
\begin{array}{lll}
-\Delta W=1 \,, &\mbox{in }V_\delta\,,\\
W= 0\,, &\mbox{on }\partial \Omega \,.\\
\end{array}\right.
\]
By comparison, we have that $\hat{w}\le W$ in $V_\delta$, which means $w(x)\le c_0\delta W(x)$\,. Moreover, we know that the function $W$ satisfies $W(x)\le c_1\,d(x)\le c_2 \Phi_1(x)$ since we know that $\Phi_1\ge c\,d(x)$. Summing up we have proved that
\[
w(x)\le c_0\delta W(x)\le c_1\delta \,d(x)\le c_2\delta \Phi_1(x)\,,
\]
recalling that $w=(1-\varepsilon)\Phi_1-U$, the above inequality gives
\[
\big[(1-\varepsilon)-c_2\delta\big]\Phi_1\le  U
\]
Now putting $c_2\delta\le \varepsilon$ we get the result, when $\delta$ is small enough.

\noindent$\bullet~$\textsc{Conclusion. }The two above steps imply that given $\ve>0$ there exist a $\delta>0$ and $p_\ve>1$ such that the above two steps hold and
\[
(1-2\ve)\Phi_1\le U_p\le (1+\ve)\Phi_1\qquad\mbox{in }\Omega\,\mbox{.\qed}
\]

Combining Proposition \ref{prop.47} with Proposition \ref{prop.k0.k1} we get the full statement of Theorem \ref{Elliptic.THM}.\qed
\subsection{Additional bounds on $\lambda_p$}
\noindent We shall also prove suitable lower bounds for $\lambda_p$, both for the sake of completeness and because they will be used in Section \ref{ssec.state}. These bounds are easier to obtain than the upper bounds.

\noindent(i) Using $U_p$  as test function, we obtain the global energy equality $\lambda_p\|U_p\|_{p+1}^{p+1}=\|\nabla U_p\|_2^2$, that combined with the Sobolev inequality
\[
\|f\|_{p+1}\le |\Omega|^{\frac{1}{p+1}-\frac{1}{2^*}}\|f\|_{2^*}\le |\Omega|^{\frac{1}{p+1}-\frac{1}{2^*}}\mathcal{S}_2\|\nabla f\|_2
\]
gives, recalling that we have chosen $\lambda_p$ in such a way that $\|U_p\|_{p+1}=1$,
\[
\frac{1}{|\Omega|^{\frac{2}{p+1}-\frac{2}{2^*}}}=\frac{\|U_p\|_{p+1}^{2}}{|\Omega|^{\frac{2}{p+1}-\frac{2}{2^*}}}\le \left[\int_{\Omega} U_p^{2^*}\dx\right]^{\frac{2}{2^*}}
     \le \mathcal{S}_2^2\|\nabla U_p\|_2^2= \mathcal{S}_2^2\lambda_p\|U_p\|_{p+1}^{p+1}= \mathcal{S}_2^2\lambda_p\,.\\
\]
We can rewrite the lower bound as follows
\begin{equation}\label{lower.lambda.p.0}
\frac{1}{\mathcal{S}_2^2\,|\Omega|^{\frac{2}{p+1}-\frac{2}{2^*}}}\le \lambda_p\qquad\mbox{and for } p\to 1\qquad \frac{1}{\mathcal{S}_2^2\,|\Omega|^{1-\frac{2}{2^*}}}\le \lambda_1\,.
\end{equation}
\noindent(ii) Other lower bounds can be obtained by combining H\"older, Poincar\'e and Sobolev inequalities:
\[
\|U_p\|_{p+1}^2\le \|U_p\|_{2^*}^{2\vartheta}\|U_p\|_2^{2(1-\vartheta)}
    \le \frac{\left(\lambda_1\mathcal{S}_2^2\right)^{\vartheta}}{\lambda_1}
        \|\nabla U_p\|_2^2\qquad\mbox{with}\qquad\vartheta=\frac{d(p-1)}{2(p+1)}
\]
which gives
\begin{equation}\label{lower.ii}
\lambda_p= \int_\Omega|\nabla U_p|^2\dx \ge \frac{\lambda_1}{\left(\lambda_1\mathcal{S}_2^2\right)^{\vartheta}}\|U_p\|_{p+1}^2 =\lambda_1\left(\lambda_1\mathcal{S}_2^2\right)^{-\frac{d(p-1)}{2(p+1)}}
\end{equation}
since we have chosen $\lambda_p$ in such a way that $\|U_p\|_{p+1}=1$.

\noindent\textbf{The case of variational solutions. }Other estimates for $\lambda_p$ can be easily obtained in the case in which solutions are minima of a suitable functional, this happens for instance in the case of domains $\Omega$ for which the solution is unique, hence they are minima, since a solution which is a minima always exists as a consequence of Kondrachov's compactness theorem.

When the solution of the Elliptic problem \ref{ell.prob.p} are minima of a suitable functional, namely
when we consider the homogeneous functional
\begin{equation*}
J_p[u]=\frac{\int_\Omega|\nabla u|^2\dx}{\left(\int_\Omega u^{p+1}\dx\right)^{\frac{2}{p+1}}}
\end{equation*}
defined on $W_0^{1,2}(\Omega)$, and we seek for its minimum under the restriction $\|u\|_{p+1}=1$, we can define
\begin{equation*}
\lambda_p=\inf_{u\in X_p}J_p[u]=\inf_{u\in X_p}\int_\Omega|\nabla u|^2\dx\qquad\mbox{where}\qquad X_p=\left\{u\in W_0^{1,2}(\Omega)\;\big|\;\|u\|_{p+1}=1  \right\}\,.
\end{equation*}
Let $U_p\in X_p$ be a solution to the elliptic problem \ref{ell.prob.p} with $\lambda_p$ defined as above. Estimates in this case are simpler and hold for any $1\le p< p_s$.
\begin{prop}\label{lambda.prop}
Under the above assumptions, if $U_p$ is a minimum for the functional $J_p$ on the set $X_p$, then  it is a positive weak (hence classical) solution to the elliptic Problem {\rm \ref{ell.prob.p}}. Moreover the following estimates hold
\begin{equation}\label{lambda.1.p}
\left(\mathcal{S}_2\lambda_1\right)^{-\frac{d(p-1)}{2(p+1)}}\le \frac{\lambda_p}{\lambda_1}=\frac{\inf\limits_{u\in X_p}J_p[u]}{\inf\limits_{u\in X_1}J_1[u]}\le |\Omega|^{\frac{p-1}{p+1}}
\end{equation}
where $\lambda_1$ is the first eigenvalue of the Dirichlet Laplacian on $\Omega$, and $\mathcal{S}_2$ is the constant on the Sobolev imbedding from $W^{1,2}_0(\Omega)$.  As a consequence, $\lambda_p\to\lambda_1$ as $p\to 1^+$.
\end{prop}
\noindent {\sl Proof.~}It is a standard fact in calculus of variations to see that a minimum of $J_p$ is a weak solution to the elliptic problem under consideration. We can now prove the upper estimate:
\[
\lambda_p=\inf_{u\in X_p}\int_\Omega|\nabla u|^2\dx = \inf_{u\in W^{1,2}_0(\Omega)} \frac{\int_\Omega|\nabla u|^2\dx}{\left(\int_\Omega u^{p+1}\dx\right)^{\frac{2}{p+1}}}
\le \frac{\int_\Omega|\nabla \Phi_1|^2\dx}{\left(\int_\Omega \Phi_1^{p+1}\dx\right)^{\frac{2}{p+1}}}\le \lambda_1\, |\Omega|^{\frac{p-1}{p+1}}
\]
if we moreover assume $\|\Phi_1\|_2=1$ (not restrictive). We have just used the fact that $\Delta\Phi_1=\lambda_1\Phi_1$ together with H\"older inequality $\|\Phi_1\|_2^2\le |\Omega|^{\frac{p-1}{p+1}}\|\Phi_1\|_{p+1}^2$. The lower estimates are exactly the same as \eqref{lower.ii} and we do not repeat the proof here.\qed

\section{Convergence with rates for $m$ near one}\label{sec.Stabnear1}

The idea of this section is simple, but the technical details lengthy. Since the lower bound for the decay rates in the linear case $m=1$ is just  $\lambda_0 {\bf c}=\lambda_2-\lambda_1>0$, it must be also positive for $m$ near 1 by continuity. Putting the details into this program is not so easy and we present below the part that we have been able to prove. The section starts by proving a suitable Poincar\'e inequality. It continues by estimating the constant ${\bf c}$ that enters the elliptic problem and proving that it tends to $\lambda_1$ as it should. These two ingredients allow to state and prove our main results about convergence with rate.

\subsection{Weighted Poincar\'e Inequality}
We recall that putting $p=1/m$ and $S_m=u_p^p$ we get a solution $S_m$ to
\[
\left\{
\begin{array}{lll}
-\Delta S_m^m=\lambda_m\,S_m&\mbox{in }\Omega\\
S_m>0&\mbox{in }\Omega\\
S_m=0&\mbox{on }\partial\Omega\\
\end{array}
\right.
\]
for which we know, by Theorem \ref{Elliptic.THM}, that $1=\|S_m\|_{m+1}=\|\Phi_1\|_2$ and, as $m\to 1^-$, $\lambda_m\to \lambda_1$, and $\widetilde{k}_0(1/m)\Phi_1\le S_m^m\le \widetilde{k}_0(1/m)\Phi_1 $, with $\lim_{m\uparrow 1}\widetilde{k}_i(1/m)=1$.  Hence setting $k_i(m):=\widetilde{k}_i(1/m)$ we have:

\medskip

\noindent $\mathbf{(H_m)}$ For any $m_s = (d-2)/(d+2)<m\le 1$ there exists constants, $k_i(m)$ with $\lim_{m\uparrow1}k_i(m)=1$, such that the stationary solutions $S_m(x)$ satisfy the bound
\begin{equation}\label{H1.p.bis}
k_0(m)\,\Phi_1(x)\le S_m^m(x) \le k_1(m)\,\Phi_1(x)\qquad\mbox{for any }x\in\overline{\Omega}\,.
\end{equation}

\begin{thm}[Weighted Poincar\'e Inequality]\label{Weighted.Poincare} Let $f\in W_0^{1,2}(\Omega)$ and $g=f/\Phi_1$. Let $S_m$ be a weight satisfying $\mathbf{(H_m)}$. Then the following inequality holds
\begin{equation}\label{Weight.Poincare}
\frac{\Lambda\,k_0(m)^2}{k_1(m)^2\|S_m\|_\infty^{1-m}}\int_\Omega \left|g-\overline{g}\right|^2 S_m^{1+m} \dx
    \le \int_\Omega \left|\nabla g \right|^2 S_m^{2m}\dx
\end{equation}
where $\Lambda=\lambda_2-\lambda_1>0$ is the optimal constant in the intrinsic Poincar\'e inequality \eqref{Poincare.w} and
\[
\overline{g}=\frac{\int_\Omega g S_m^{1+m} \dx}{\int_\Omega S_m^{1+m} \dx}\,.
\]
\end{thm}
\noindent {\sl Proof.~}Notice first that, by $\mathbf{(H_m)}$, $\Phi_1^2(x)\le S_m^{2m}/k_0(m)^2$.
As a consequence,
\begin{equation}\label{up}
\int_\Omega \left|\nabla g \right|^2\Phi_1^2\dx
    \le \frac{1}{k_0(m)^2} \int_\Omega \left|\nabla g \right|^2 S_m^{2m}\dx\,.
\end{equation}
Moreover
\begin{equation}\label{lower}
\Phi_1^2(x)\ge \frac{S_m^{2m}}{k_1(m)^2}\ge \frac{S_m^{1+m}}{k_1(m)^2\|S_m\|_\infty^{1-m}}\,,
\end{equation}
where we used again $\mathbf{(H_m)}$ and the fact that $S_m^{2m}(x)\ge \|S_m\|_\infty^{m-1}S^{1+m}\,$, valid since $m<1$. Therefore
\begin{equation}\label{low}\begin{split}
\int_\Omega \left|g-g_{\Phi_1}\right|^2 \Phi_1^2 \dx
    &\ge \frac{1}{k_1(m)^2\|S_m\|_\infty^{1-m}}\int_\Omega \left|g-g_{\Phi_1}\right|^2 S_m^{1+m} \dx\\
    &\ge \frac{1}{k_1(m)^2\|S_m\|_\infty^{1-m}}\int_\Omega \left|g-\overline{g}\right|^2 S_m^{1+m} \dx
\end{split}
\end{equation}
where the last inequality follows  by Lemma \ref{lemma.mean.opt}.

\noindent Putting together the latter inequalities yields \eqref{Weight.Poincare}
\begin{equation*}\begin{split}
\frac{\Lambda}{k_1(m)^2\|S_m\|_\infty^{1-m}}\int_\Omega \left|g-\overline{g}\right|^2 S_m^{1+m} \dx
    &\le \Lambda\int_\Omega \left|g-g_{\Phi_1}\right|^2 \Phi_1^2 \dx\\
    &\le \int_\Omega \left|\nabla g \right|^2\Phi_1^2\dx
     \le \frac{1}{k_0(m)^2} \int_\Omega \left|\nabla g \right|^2 S_m^{2m}\dx\mbox{\,.\qed}
\end{split}\end{equation*}

We now recall the following well--known result which has been used in the above proof. Hereafter, $f_\mu:=\mu(X)^{-1}\int_Xf\,{\rm d}\mu$ where $\mu$ is any nonnegative bounded measure.
\begin{lem}\label{lemma.mean.opt}
Let $f\in \LL^2(X,\rd\mu)$, with $\mu(X)<\infty$. Then we have
\begin{equation}\label{mean-const.2}
\left\|f-{f}_\mu\right\|_{\LL^2(X,\rd\mu)}\le \left\|f- c\right\|_{\LL^2(X,\rd\mu)}\,,
    \qquad\mbox{for all $c\in \RR$\,.}
\end{equation}
\end{lem}
\noindent {\sl Proof.~}
By expanding the square
\[\begin{split}
\left\|f-c\right\|_{\LL^2(X,\rd\mu)}^2
    &= \int_X |f-c|^2 \rd\mu
     =\mu(X)\int_X \big[f^2 -2cf +c^2\big] \frac{\rd\mu}{\mu(X)}\\
    &=\mu(X)\left[ \int_X f^2 \frac{\rd\mu}{\mu(X)}
        -2c \int_X  f \frac{\rd\mu}{\mu(X)} + c^2 \right]
             =\mu(X)\left[\int_X f^2 \frac{\rd\mu}{\mu(X)}
        -2c {f}_\mu + c^2 \right]\\
    &\ge\mu(X)\left[\frac{\left\|f\right\|_{\LL^2(X,\rd\mu)}^2}{\mu(X)}-\big({f}_\mu\big)^2\right]
     =\left\|f\right\|_{\LL^2(X,\rd\mu)}^2- \mu(X)\big({f}_\mu\big)^2
     =\left\|f-{f}_\mu\right\|_{\LL^2(X,\rd\mu)}^2\mbox{\,.\qed}
\end{split}
\]

The version we will really need  is a variation in the use of different weights.
Let  $S_{c,m}$ be a positive solution of $-\Delta S^m = {\bf c} S$ on $\Omega$, vanishing at the boundary, obtained
as asymptotic profile for an evolution with fixed initial data $u_0$ and variable $m<1$.

\begin{thm}\label{Weighted.Poincare.c} Let $f\in W_0^{1,2}(\Omega)$ and $g=f/\Phi_1$. Let $S_{c,m}$ be as above. Then the following inequality holds
\begin{equation}\label{Weight.Poincare.c}
\frac{\Lambda\,k_0(m)^2}{k_1(m)^2\|S_{c,m}\|_\infty^{1-m}}\int_\Omega \left|g-\overline{g}\right|^2 S_{c,m}^{1+m} \dx
    \le \int_\Omega \left|\nabla g \right|^2 S_{c,m}^{2m}\dx
\end{equation}
where $\Lambda=\lambda_2-\lambda_1>0$ is the optimal constant in the intrinsic Poincar\'e inequality \eqref{Poincare.w} and
\[
\overline{g}=\frac{\int_\Omega g S_{c,m}^{1+m} \dx}{\int_\Omega S_{c,m}^{1+m} \dx}\,,
\]
and we know that  $\lim_{m\uparrow1}k_1(m)/k_0(m)=1$.
\end{thm}

\noindent {\sl Proof.~}We use the preceding result and the following observation: it easy to check the relations between $S_m$ and $S_{c,m}$:
\[
S_{c,m}=\left(\frac{\lambda_m}{{\bf c}}\right)^{\frac{1}{1-m}} S_m\qquad\mbox{and}\qquad
\|S_{c,m}\|_q^{1-m} =\frac{\lambda_m}{{\bf c}}\|S_m\|_q^{1-m}
\]
We now plug the above equalities in the weighted Poincar\'e \eqref{Weight.Poincare} to get
\[
\frac{\Lambda\,k_0(m)^2}{k_1(m)^2 \frac{{\bf c}}{\lambda_m}\|S_{c,m}\|_\infty^{1-m}}
    \int_\Omega \left|g-\overline{g}\right|^2  \left(\frac{{\bf c}}{\lambda_m}\right)^{\frac{1+m}{1-m}} S_{c,m}^{1+m}\dx
        \le \int_\Omega \left|\nabla g \right|^2 \left(\frac{{\bf c}}{\lambda_m}\right)^{\frac{2m}{1-m}} S_{c,m}^{2m}\dx
\]
that is exactly \eqref{Weight.Poincare.c} since the factors ${\bf c}/\lambda_m$ simplify.\qed

\subsection{Estimating the extinction time and the constant $\bf c$}\label{estimates.T.1-m}
We need to estimate the extinction time $T=T(m,d,u_0)$ from above and from below, to obtain bounds on the constant ${\bf c}=1/(1-m)T$ appearing in the rescaled equation \eqref{FDE.Problem.Domain.Fine.Log} and the elliptic equation \eqref{FDE.Elliptic.Problem}.

\begin{prop}\label{prop.C}Let $0<m<1$ and $u$ be the solution to the original Problem ~{\rm \ref{FDE.Problem.Domain.Fine}} corresponding to an initial datum $u_0\in \LL^r(\Omega)$ with $r>1$ and $r\ge r_c=d(1-m)/2$. Then its extinction time $T=T(m,d,u_0)$ satisfies the bounds
\begin{equation}\label{bounds.T}
\frac{1}{\lambda_1}\frac{\left[\int_\Omega u_0(x) \Phi_1(x)\dx\right]^{1-m}}{\left[\int_\Omega\Phi_1(x)\dx\right]^{1-m}}
\le (1-m) T \le \frac{(r+m-1)^2}{4m(r-1)}\frac{\left(\lambda_1\mathcal{S}_2^2\right)^{\frac{d(1-m)}{4r}}}{\lambda_1} \|u_0\|_r^{1-m}\,.
\end{equation}
Taking $r=1+m$, which amounts to ask $m>m_s=(d+2)/(d-2)$, we get
\begin{equation}
\frac{1}{\lambda_1}\frac{\left[\int_\Omega u_0(x) \Phi_1(x)\dx\right]^{1-m}}{\left[\int_\Omega\Phi_1(x)\dx\right]^{1-m}}
\le (1-m) T \le \frac{\left(\lambda_1\mathcal{S}_2^2\right)^{\frac{d(1-m)}{4(1+m)}}}{\lambda_1} \|u_0\|_{1+m}^{1-m}\,.
\end{equation}
\end{prop}

\begin{cor}\label{cor55} If $u_0\in L^{m+1}(\Omega)$ and $m>m_s$  we have
\begin{equation}\label{limit.T(1-m)}
\lim_{m\to 1^-}(1-m)T(m,d,u_0)=\frac{1}{\lambda_1}\,.
\end{equation}
hence ${\bf c}\to \lambda_1$ as  $m\to 1$. Moreover, $({\bf c}/\lambda_1)^{1/{(1-m)}}=O(1)$ as $m\to 1$.
\end{cor}

\noindent {\sl Proof of the Proposition.~}We begin with the lower bound. We take $\Phi_1$ as test function and consider the solution $u$ to the original Problem \ref{FDE.Problem.Domain.Fine}
\begin{equation*}
\left\{\begin{array}{lll}
u_\tau=\Delta (u^m) & ~ {\rm in}~ (0,+\infty)\times\Omega\\
u(0,x)=u_0(x) & ~{\rm in}~ \Omega \\
u(\tau,x)=0 & ~{\rm for}~  \tau >0 ~{\rm and}~ x\in\partial\Omega\\
\end{array}\right.
\end{equation*}
and derive the integral
\begin{equation*}
\begin{split}
\left|\frac{\rd}{\rd\tau}\int_\Omega u(\tau,x)\Phi_1(x)\dx\right|
&=\left|\int_\Omega\left( \Delta u^m(\tau,x)\right)\Phi_1(x)\dx\right|
 =\left|\int_\Omega u^m(\tau,x) \Delta\Phi_1(x)\dx\right|
 =\left|\int_\Omega u^m(\tau,x) \lambda_1\Phi_1(x)\dx\right|\\
&= \lambda_1 \int_\Omega u^m(\tau,x) \Phi_1(x)\dx
\le \lambda_1 \left[\int_\Omega u(\tau,x) \Phi_1(x)\dx\right]^m \left[\int_\Omega\Phi_1(x)\dx\right]^{1-m} \\
\end{split}
\end{equation*}
where we have integrated by parts since both $u$ and $\Phi_1$ are zero at $\partial\Omega$, and we recall that $\lambda_1\Phi_1=-\Delta\Phi_1$. Integrating the differential inequality gives
\[
\left[\int_\Omega u(t,x) \Phi_1(x)\dx\right]^{1-m}
\le \left[\int_\Omega u(s,x) \Phi_1(x)\dx\right]^{1-m}+\lambda_1(1-m)\left[\int_\Omega\Phi_1(x)\dx\right]^{1-m}|t-s|
\]
for any $0\le s, t\le T$. Letting $s=T$ and $t=0$ gives the lower bound \eqref{bounds.T}.

\noindent The upper bound follows  by using H\"older, Sobolev and Poincar\'e inequalities in the form
\begin{equation*}
\|f\|_{q}^2\le \|f\|_{2^*}^{2\vartheta}\|f\|_2^{2(1-\vartheta)}
    \le \frac{\left(\lambda_1\mathcal{S}_2^2\right)^{\vartheta}}{\lambda_1}
        \|\nabla f\|_2^2\qquad\mbox{with}\qquad\vartheta=\frac{d(q-2)}{2q}\qquad\mbox{and}\qquad 2\le q \le2^*
\end{equation*}
We are going to use the above inequality for $q=(2r)/(r+m-1)$ and $f=u^{\frac{r+m-1}{2}}$, noticing that $q\in [2,2^*]$ if and only if $r\ge r_c=d(1-m)/2$, which is
\[
\frac{\lambda_1}{\left(\lambda_1\mathcal{S}_2^2\right)^{\frac{d(1-m)}{4r}}}\left\|u\right\|_r^{r+m-1}\le
        \left\|\nabla u^{\frac{r+m-1}{2}}\right\|_2^2
\]
Differentiation of the $L^{r}-$norm gives when $r>1$ and $r\ge r_c=d(1-m)/2$
\begin{equation}\label{Lp.Global.Time.Deriv}
\begin{split}
\frac{\rd}{\dt}\|u(t)\|_r^r
    &=-\frac{4mr(r-1)}{(r+m-1)^2}\int\left|\nabla u(t)^{\frac{r+m-1}{2}}\right|^2\dx
        \le -\frac{4mr(r-1)}{(r+m-1)^2}\frac{\lambda_1}
            {\left(\lambda_1\mathcal{S}_2^2\right)^{\frac{d(1-m)}{4r}}}\left\|u(t)\right\|_r^{r\left(1-\frac{1-m}{r}\right)}\\
\end{split}
\end{equation}
which is a closed differential inequality of the form $Y'\le -K Y^{1-\varepsilon}$, which integrated on $[s,t]$ gives
\[
\|u(t)\|_r^{1-m}-\|u(s)\|_r^{1-m}\le - \frac{4mr(r-1)}{(r+m-1)^2}\frac{\lambda_1}
            {\left(\lambda_1\mathcal{S}_2^2\right)^{\frac{d(1-m)}{4r}}}\frac{1-m}{r}(t-s)
\]
which gives the upper bound, just by letting $s=0$ and $t=T$.\qed

\subsection{Statement of the main convergence result}\label{ssec.state}
The next step consists in showing that these inequalities allow us to apply the decay results of Subsection \ref{ssec.frc}.
For that we need to estimate in a clear way the constant before the left integral in \eqref{Weight.Poincare.c} and obtain a lower bound with a constant $K$ independent of the particular solution $S$. This is possible thanks to the following estimate
\[
\|S_{c,m}\|_\infty^{1-m}=\frac{\lambda_m}{{\bf c}}\|S_m\|_\infty^{1-m}
\ge \frac{\lambda_m}{{\bf c}}\frac{\|S_m\|_{m+1}^{1-m}}{|\Omega|^{\frac{1-m}{1+m}}}
=\frac{\lambda_m}{{\bf c}}\frac{1}{|\Omega|^{\frac{1-m}{1+m}}}\ge\frac{\lambda_1}{{\bf c}}
\,\left[\left(\mathcal{S}_2\lambda_1\right)^{\frac{d}{2}}\,|\Omega|\right]^{\frac{m-1}{1+m}}
\]
which is true since $\|S_m\|_{1+m}^{1-m}=1$ by construction and since the lower bound \eqref{lower.ii}, rewritten for $p=1/m$, reads
\begin{equation*}
\frac{\lambda_m}{\lambda_1}\ge \left(\mathcal{S}_2\lambda_1\right)^{-\frac{d(1-m)}{2(1+m)}}\,.
\end{equation*}
This means that the Weighted Poincar\'e Inequality mentioned in Subsection \ref{ssec.frc} holds in the form
\begin{equation}\label{ewpi2.1}
K{\bf c} \dint_\Omega S_c^{m+1}|g-\overline g|^2\,dx\le \dint _\Omega S_c^{2m}|\nabla g|^2\,dx\,,
\end{equation}
with
\[
K=K(m)=\frac{(\lambda_2-\lambda_1)\,k_0(m)^2 }{\lambda_1 k_1(m)^2}\left[\left(\mathcal{S}_2\lambda_1\right)^{\frac{d}{2}}\,|\Omega|\right]^{\frac{1-m}{1+m}}.
\]
At this moment, we see that the necessary condition to obtain decay is then
$$
F(m) := m K(m)-2(1-m)>0,
$$
in other words,  $m>2/(2+K(m))$. But since  $F(1)=(\lambda_2-\lambda_1)/\lambda_1>0$, it follows that there is an $m_\sharp < 1$ such that for all $m_\sharp<m<1$ we have that $F(m)>0$. Note that $m_{\sharp}$ changes with the geometry of the domain. It may be objected that $m_{\sharp}$ is given in a very implicit way. However, given an estimate of the form  $k_1(m)/k_0(m)\le C$ when $m_c<m<1$, then $1-m_\sharp$ can be explicitly estimated from below in terms of $C$, $\lambda_1$, $\lambda_2$, $\mathcal{S}_2$ and $|\Omega|$. We shall provide suitable explicit bounds for $k_i(m)$ in a forthcoming paper \cite{BGV-Elliptic}, see also Theorem \ref{prop.k0.k1.Appendix}.

We can now state the rescaled version of the asymptotic convergence result. The rate will involve the expression
\begin{equation}\label{rate.condition.intro}
\gamma_0(m)=\frac{1}{(1-m)T}\left[m\left(\frac{\lambda_2}{\lambda_1}-1\right)\frac{k_0(m)^2}{k_1(m)^2}
    \left[\left(\mathcal{S}_2\lambda_1\right)^{\frac{d}{2}}\,|\Omega|\right]^{\frac{1-m}{1+m}}
        -2(1-m)\right]>0
\end{equation}
where $m_\sharp<m<1$, the constants $k_i(m)\to 1$ as $m\to 1$, cf. Theorem \ref{Elliptic.THM} (rewritten there with $p=1/m$). We recall that by Proposition \ref{prop.C} the quantity $(1-m)T$ appearing in \eqref{rate.condition.intro} can be explicitly estimated from above and below in terms of $u_0$ and moreover, by Corollary \ref{cor55}, we have that $\lim\limits_{m\to 1^-}1/(1-m)T=\lambda_1$.

\begin{thm}[Rates of convergence, rescaled version]\label{thm.main.decay.rescaled}
Let $\max\{m_{\sharp},m_c\}< m <1$. Let $v$ be the rescaled solution corresponding to an initial datum $u_0$ as in Theorem {\rm \ref{Thm.Main.Intro}}, which converge to its unique stationary profile $S$. Let $\gamma<\gamma_0$ then for all $t>t_0$, with $t_0$ sufficiently large, we have the following entropy decay formula:
\begin{equation}\label{rate.entropy.rescaled}
\mathcal{E}[\theta(t)]\le \ee^{-\gamma(t-t_0)}\mathcal{E}[\theta(t_0)]\,.
\end{equation}
In other words,  the weighted $\LL^2$-norm decays with rate $\gamma$, more precisely there exists constants $\kappa_i>0$ and a time $t_0>0$ such that
\begin{equation}\label{rate.intro.rescaled}
\int_\Omega \left|v(t,x)-S(x)\right|^2 S(x)^{m-1}\dx=\int_\Omega\left|\frac{v(t,x)}{S(x)}-1\right|^2 S(x)^{1+m}\dx
    \le \kappa_0\,\ee^{-\gamma(t-t_0)}\,.
\end{equation}
Moreover for all $q\in(0,\infty]$ \begin{equation}\label{rate.intro.infty.rescaled}
\|v(t,\cdot)-S(\cdot)\|_q \le\,\kappa_1\,\ee^{-\frac{\gamma}{2}(t-t_0)}
\end{equation}
where the constant $\kappa_1$  depends on $m,d$ and $u_0$ and on the uniform bounds on the $C^\alpha-$norm of $u(t_0)$.
\end{thm}
\noindent\textbf{Remark. }The constant $\gamma_0$ satisfies
\[
\lim_{m\to 1^+}\gamma_0(m)=\lambda_2-\lambda_1\,,
\]
as a consequence of Corollary \ref{cor55}.
\begin{thm}[Rates of convergence, original variables]\label{thm.main.decay}
Let $\max\{m_{\sharp},m_c\}< m <1$. Let $u$ be the solution to Problem \ref{FDE.Problem.Domain.Fine}, let $T=T(m,d,u_0)$ be its extinction time and let $\mathcal{U}_T$ be as in Theorem {\rm \ref{Thm.Main.Intro}}, so that $u(\tau)/\mathcal{U}_T(\tau)\to 1$ as $\tau\to T$. Then, for any
\[
{\overline \gamma}< \frac{1}{(1-m)}\left[m\left(\frac{\lambda_2}{\lambda_1}-1\right)\frac{k_0(m)^2}{k_1(m)^2}
\left[\left(\mathcal{S}_2\lambda_1\right)^{\frac{d}{2}}\,|\Omega|\right]^{\frac{1-m}{1+m}}
        -2(1-m)\right]
\]
there exists a constant $\kappa>0$ such that
\begin{equation}\label{rate.intro}
\left\|\frac{u(\tau,\cdot)}{\mathcal{U}(\tau,\cdot)}-1\right\|_{\LL^2(\Omega, S^{1+m})}^2
\le \kappa_0 \,\left(\frac{T-\tau}{T}\right)^{{\overline\gamma}}
\end{equation}
or equivalently
\begin{equation}\label{rate.intro.equiv}
\int_\Omega \left|u(\tau,x)-\mathcal{U}(\tau,x)\right|^2 S^{m-1}\dx\le \kappa_0 \,\left(\frac{T-\tau}{T}\right)^{\frac{2}{1-m}+{\overline\gamma}}
\end{equation}
for all $t_0\le \tau \le T$, where $\kappa_0$ depends on $m,d$ and $u_0$.  Moreover we have that for all $q\in(0,\infty]$
\begin{equation}\label{rate.infty}
\|u(\tau,x)-\mathcal{U}(\tau,x)\|_q
    \le \kappa_1 \,\left(\frac{T-\tau}{T}\right)^{\frac{2}{1-m}+{\overline\gamma}}
\end{equation}
where the constant $\kappa_1$  depends on $m,d$ and $u_0$ and on the uniform bounds on the $C^\alpha-$norm of $u(t_0)$.
\end{thm}

We comment that the weighted convergence of \eqref{rate.intro.equiv} is somehow stronger than the non-weighted $\LL^p-$norm convergence, since the weight $S^{m-1}$ is singular at the boundary.

\medskip

\noindent The main result for the Porous Medium Equation reads:
\begin{thm}\label{PME.Rates} Let $m>1$, let $v$ be a the rescaled solution as in Subsection \ref{Sect.PME} that converges to its unique stationary state $S$, and let $\theta=v/S$. Then, for all $0<\beta<2+\frac{Km}{m-1}$ there exists a time $t_1$ depending on $m,d,\beta$ and on the constant $K>0$ of the GWPI, such that the entropy decays as
\begin{equation}
\mathcal{E}[\theta(t)]\le \mathcal{E}[\theta(t_1)]\,\ee^{-\beta(t-t_1)}\qquad\mbox{for all }t\ge t_1.
\end{equation}
Moreover for all $q\in(0,\infty]$
\begin{equation*}
\|v(t,\cdot)-S(\cdot)\|_{\LL^q(\Omega)} \le\,\kappa_1\,\ee^{-(t-t_0)}
\end{equation*}
for all $t>t_1\gg 1$, where the constant $\kappa_2$  depends on $m,d$ and $u_0$ and on the uniform bounds on the $C^\alpha-$norm of $u(t_0)$. In original variables we obtain that for all $q\in(0,\infty]$
\[
\left\|u(\tau,\cdot)-\mathcal{U}(\tau,\cdot)\right\|_{\LL^q(\Omega)}\le \frac{\kappa_2}{(1+\tau)^{1+\frac{1}{m-1}}}\,.
\]
\end{thm}

\medskip

In order to conclude the proof of the above theorems, we need an interpolation Lemma due to Gagliardo \cite{Ga}, cf. also  Nirenberg, \cite[p. 126]{MR0109940}.
\begin{lem}\label{GN} Let $\lambda$, $\mu$ and $\nu$ be such that $-\infty<\lambda \le \mu
\le \nu <\infty$. Then there exists a positive constant $\mathcal
C_{\lambda,\mu,\nu}$ independent of $f$ such that
\begin{equation}\label{eq:interpolation}
\|f\|_{1/\mu}^{\nu-\lambda} \le \mathcal
C_{\lambda,\mu,\nu}\|f\|_{1/\lambda}^{\nu-\mu} \;
\|f\|_{1/\nu}^{\mu-\lambda}\quad\forall\;f\in\mathcal C(\RR^d)\;,
\end{equation}
where $\|\cdot\|_{1/\sigma}$ stands for the following quantities:\begin{itemize}
\item[(i)] If $\sigma>0$, then
$\|f\|_{1/\sigma}=\left(\int_{\RR^d}|f|^{1/\sigma}\dx\right)^\sigma$.

\item[(ii)] If $\sigma<0$, let $k$ be the integer part of
$(-\sigma d)$ and $\alpha=|\sigma|d-k$ be the fractional
(positive) part of $\sigma$. Using the standard multi-index
notation, where $|\eta|=\eta_1+\ldots+\eta_d$ is the length
of the multi-index
$\eta=(\eta_1,\ldots\eta_d)\in\mathbb{Z}^d$, we define
\begin{equation*}\label{def.C^k}
\|f\|_{1/\sigma}=\left\{\begin{array}{lll} \displaystyle\max_{|\eta|=k}\; \big|\partial^\eta f\big|_\alpha=\displaystyle\max_{|\eta|=k}\; \sup_{x,y\in\RR^d}\;\dfrac{\big|\partial^\eta f(x)-\partial^\eta f(y)\big|}{|x-y|^\alpha}=|f\|_{C^\alpha(\RR^d)}& \mbox{if~}\alpha>0\;,\\[5mm]
\displaystyle\max_{|\eta|=k}\;\displaystyle\sup_{z\in\RR^d}\big|\partial^\eta f(z)\big|:=\|f\|_{C^k(\RR^d)}& \mbox{if~}\alpha=0\;.
\end{array}\right.
\end{equation*}
As a special case, we observe that $\|f\|_{-d/j}=\|f\|_{C^{j}(\RR^d)}$.

\item[(iii)] By convention, we note $\|f\|_{1/0}=\sup_{z\in\RR^d}|f(z)|=\|f\|_{C^0(\RR^d)}=\|f\|_{\infty}$.
\end{itemize}
\end{lem}

Next we need a regularity result that helps us to compare the $C^{\alpha}$-norm with the $\LL^\infty$-norm, and we combine it with the above interpolation in order to obtain the same rate of decay for all $\LL^p$-norms.

\begin{lem}\label{Lemma.Calpha} Let $m>0$. Let $v$ be the rescaled solution  to equation \eqref{FDE.Problem.Domain.Fine.Log} corresponding to an initial datum $u_0$ as in Theorem {\rm \ref{Thm.Main.Intro}}. There exists $t_0\geq 0$, $\alpha \in (0,1)$ and a constant $B>0$ such that $v(t,x)-S(x)$ is in $C^\alpha(\Omega)$ and
\begin{equation}\label{eq:mathcalH}
\|v(t,x)-S(x)\|_{C^\alpha(\Omega)}\le\,B\,\|v(t,x)-S(x)\|_{\LL^\infty(\Omega)}\quad\forall\; t\geq t_0\;.
\end{equation}
\end{lem}
\noindent {\sl Proof.~}Let $h(t,x):=v(t,x)-S(x)$. Since both $v$ and $S$ are solutions to equation \eqref{FDE.Problem.Domain.Fine.Log}, $h$ solves
\begin{equation*}\label{Eq.Difference}
h_t = (v-S)_t = v_t= \Delta (v^m)+{\bf c}v =\Delta \big((h+S)^m\big)+{\bf c}(h+S)\,.
\end{equation*}
By Theorem~\ref{Thm.Main.Intro}, we know that for some $t_0\geq 0$, for any $t\geq t_0$, $\|h(t)\|_\infty$ can be taken uniformly small and $v$ is positive in $\Omega$, $v=0$ on $\partial\Omega$. The H\"older continuity now follows by the nowadays classical results of DiBenedetto et al. (cf. the book \cite{DiBbook} Chap III, Thm. 1.1 for $m\ge 1$ and Chap IV, Thm. 1.1 for $0<m<1$), and holds for a class of equations of the type $h_t=\nabla\cdot A(t,x,h,\nabla h) +B(t,x,h,\nabla h)$, which
satisfy standard structure conditions:
\[
\begin{split}
A(t,x,h,\nabla h)\cdot\nabla h\ge c_0|h|^{m-1}|\nabla h|^2 -\varphi_0(x,t)\\
|A(t,x,h,\nabla h)|\le c_1|h|^{m-1}|\nabla h|^2+\varphi_1(x,t)\\
|B(t,x,h,\nabla h)|\le c_2|\nabla |h|^m |^2+\varphi_2(x,t)\\
\end{split}
\]
for suitable $c_i>0$ and nonnegative $\varphi_i$. In our case we have that
\[
A(t,x,h,\nabla h)=m(h+S)^{m-1}\nabla h\qquad\mbox{and}\qquad
B(t,x,h,\nabla h)=m\nabla\cdot((h+S)^{m-1}\nabla S)+{\bf c}(h+S) 
\]
clearly satisfy the structure conditions.\qed

The same regularity estimates can be proved for the relative error, at least in the case $m>1$.

\subsection{Proof of Theorems \ref{thm.main.decay.rescaled}, \ref{thm.main.decay} and \ref{PME.Rates}}
The result of Theorem \ref{exp.decay.norm} corresponds exactly to the weighted estimate \eqref{rate.intro.rescaled} of Theorem \ref{thm.main.decay.rescaled}. It remains to prove the $\LL^\infty-$estimate \eqref{rate.infty}. To this end we combine the results of the previous lemmata. $h(t,x):=v(t,x)-S(x)$. Lemma \ref{GN} gives
\[
\|h\|_\infty\le A\,\|h\|_{C^\alpha}^{\varepsilon}\|h\|_p^{1-\varepsilon}
\]
where we take $\alpha$ the H\"older exponent of Lemma \ref{Lemma.Calpha}, and we take any $p>0$. Then $\varepsilon\in(0,1)$ and $A=A(\alpha,p,\infty)>0$ are as in Lemma \ref{GN}. The result of Lemma \ref{Lemma.Calpha} reads
\[
\|h\|_{C^\alpha}\le\,B\,\|h\|_{\infty}\quad\forall\; t\geq t_0\;.
\]
The combination of these two results gives
\begin{equation}\label{Interp+Calpha}
\|h\|_\infty\le A\,\left(B\,\|h\|_{\infty}\right)^{\varepsilon}\|h\|_p^{1-\varepsilon}
\qquad\mbox{that is}\qquad
\|v(t,\cdot)-S(\cdot)\|_\infty\le \left(A\,B^{\varepsilon}\right)^{\frac{1}{1-\varepsilon}}\|v(t,\cdot)-S(\cdot)\|_p\,.
\end{equation}
We now combine the above interpolation inequality with the exponential decay of the weighted $\LL^2-$norm of Proposition \eqref{exp.decay.norm}: there exists a constant $\kappa>0$ such that
\begin{equation*}
\int_\Omega \left|v(t,x)-S(x)\right|^2 S(x)^{m-1}\dx=\int_\Omega\big|\theta(t)\big|^2 S^{1+m}\dx
    \le \kappa\,\ee^{-\gamma(t-t_1)}\mathcal{E}[\theta(t_1)]
\end{equation*}
for all $t>t_1\gg 1$, for all $\gamma$ such that
$$
0<\gamma< \gamma_0={\bf c}\,\left[ m\left(\frac{\lambda_2}{\lambda_1}-1\right)\frac{k_0(m)^2}{k_1(m)^2}
\left[\left(\mathcal{S}_2\lambda_1\right)^{\frac{d}{2}}\,|\Omega|\right]^{\frac{1-m}{1+m}}
        -2(1-m)\right].
$$
By H\"older inequality we have
\[
\frac{\|v(t,\cdot)-S(\cdot)\|_1^2}{\|S\|_{1-m}^{1-m}}\le \int_\Omega \left|v(t,x)-S(x)\right|^2 S^{m-1}\dx
\le \kappa\,\ee^{-\gamma(t-t_1)} \mathcal{E}[\theta(t_1)]
\]
so that, combining it with\eqref{Interp+Calpha}, for $p=1$ we obtain the second inequality \ref{rate.intro.infty.rescaled} of Theorem \ref{thm.main.decay.rescaled}
\[
\|v(t,\cdot)-S(\cdot)\|_\infty^2 \le \left(A\,B^{\varepsilon}\right)^{\frac{2}{1-\varepsilon}}\|v(t,\cdot)-S(\cdot)\|_1^2
=\left(A\,B^{\varepsilon}\right)^{\frac{2}{1-\varepsilon}} \kappa\,\ee^{-\gamma(t-t_1)} \mathcal{E}[\theta(t_1)]
:=\kappa_1\ee^{-\gamma(t-t_1)}
\]
So far we have concluded the proof of Theorem \ref{thm.main.decay.rescaled} and rescaling back we have proved Theorem \ref{thm.main.decay}.

It remains to prove Theorem \ref{PME.Rates} that is the PME case $m>1$. We just remark that
\[
\|S\|_{1-m}^{1-m}=\int_\Omega \big(S^m\big)^{\frac{1-m}{m}}\dx\sim\int_\Omega \dist(x,\partial\Omega)^{-\frac{m-1}{m}}\dx
\]
the latter quantity being finite for all $m>1$, since $0<(m-1)/m<1$\,.\qed

\section{Appendix}\label{sec.appendix}

\subsection{Intrinsic Poincar\'e inequality}

We give a proof of Proposition \ref{Intrinsic.Poincare}.
Notice first that \begin{equation}\label{Poincare.w.bis}
\lambda_2 \int_\Omega f^2 \dx \le\int_\Omega\big|\nabla f\big|^2\dx
\qquad\mbox{whenever}\qquad
f_{\Phi_1}:=\frac{\int_\Omega f\,\Phi_1\dx}{\int_\Omega \Phi_1^2\dx}=0.
\end{equation}
We now apply inequality \eqref{Poincare.w.bis} to the function
\[
f=g\Phi_1-\frac{\int_\Omega g\,\Phi_1^2\dx}{\int_\Omega\Phi_1^2 \dx}\Phi_1=\big(g-{g}_{\Phi_1}\big)\Phi_1\,,
\]
for which the above orthogonality condition clearly holds. Moreover,
we have:
\begin{equation*}\label{term.1}\begin{split}
&\int_\Omega f^2 \dx
=\int_\Omega \left[g\Phi_1-\frac{\int_\Omega g\,\Phi_1^2\dx}{\int_\Omega\Phi_1^2 \dx}\Phi_1\right]^2 \dx\\
&=\int_\Omega g^2\Phi_1^2 \dx
    +\left[\frac{\int_\Omega g\,\Phi_1^2\dx}{\int_\Omega\Phi_1^2 \dx}\right]^2\int_\Omega\Phi_1^2 \dx
    -2\frac{\int_\Omega g\,\Phi_1^2\dx}{\int_\Omega\Phi_1^2 \dx}\int_\Omega g\Phi_1^2 \dx\\
&=\int_\Omega \left[g^2-\left(\frac{\int_\Omega g\,\Phi_1^2\dx}{\int_\Omega\Phi_1^2 \dx}\right)^2\right]\Phi_1^2 \dx
    =\int_\Omega \left[g^2-\left({g}_{\Phi_1}\right)^2\right] \Phi_1^2 \dx.\\
\end{split}
\end{equation*}
In addition:
\begin{equation*}\label{term.2}\begin{split}
& \int_\Omega\big|\nabla f\big|^2\dx
    =\int_\Omega\left|\nabla \left(g\Phi_1-\frac{\int_\Omega g\,\Phi_1^2\dx}
        {\int_\Omega\Phi_1^2 \dx}\Phi_1\right)\right|^2\dx\\
    &=\int_\Omega\left|\nabla g\Phi_1\right|^2\dx
        +\left[\frac{\int_\Omega g\,\Phi_1^2\dx}{\int_\Omega\Phi_1^2 \dx}\right]^2
            \int_\Omega\left|\nabla \Phi_1\right|^2\dx
    -2\frac{\int_\Omega g\,\Phi_1^2\dx}{\int_\Omega\Phi_1^2 \dx}
            \int_\Omega \nabla (g\Phi_1)\cdot\nabla\Phi_1\dx\\
    &=\int_\Omega\left|\nabla g\Phi_1\right|^2\dx
        +\left[\frac{\int_\Omega g\,\Phi_1^2\dx}{\int_\Omega\Phi_1^2 \dx}\right]^2
            \int_\Omega \lambda_1\Phi_1^2\dx
        +2\frac{\int_\Omega g\,\Phi_1^2\dx}{\int_\Omega\Phi_1^2 \dx}
            \int_\Omega g\Phi_1\Delta\Phi_1\dx\\
   &=\int_\Omega\left|\nabla g\right|^2\Phi_1^2\dx+\lambda_1\int_\Omega g^2 \Phi_1^2\dx
        +\lambda_1\left[\frac{\int_\Omega g\,\Phi_1^2\dx}{\int_\Omega\Phi_1^2 \dx}\right]^2\int_\Omega \Phi_1^2\dx
        -2\lambda_1\left[\frac{\int_\Omega g\,\Phi_1^2\dx}{\int_\Omega\Phi_1^2 \dx}\right]^2\int_\Omega \Phi_1^2\dx\\
    &=\int_\Omega\left|\nabla g\right|^2\Phi_1^2\dx+\lambda_1\int_\Omega g^2 \Phi_1^2\dx
        -\lambda_1\left[\frac{\int_\Omega g\,\Phi_1^2\dx}{\int_\Omega\Phi_1^2 \dx}\right]^2\int_\Omega \Phi_1^2\dx\\
   &=\int_\Omega\left|\nabla g\right|^2\Phi_1^2\dx
        +\lambda_1\int_\Omega \left|g-{g}_{\Phi_1}\right|^2 \Phi_1^2 \dx\\
\end{split}
\end{equation*}

Summing up we have shown that
\[\begin{split}
\lambda_2\int_\Omega \left|g-{g}_{\Phi_1}\right|^2 \Phi_1^2 \dx&=
\lambda_2\int_\Omega f^2\dx\le\int_\Omega\big|\nabla f\big|^2\dx= \int_\Omega\left|\nabla g\right|^2\Phi_1^2\dx+\lambda_1\int_\Omega \left|g-{g}_{\Phi_1}\right|^2 \Phi_1^2 \dx
\end{split}\]
which yields the desired inequality.\qed

We recall some bounds on $\lambda_2-\lambda_1$. Singer et al. \cite{Singer}, \cite{Yu}, \cite{Ling} proved that for convex domains $\Omega\subset \RR^d$ with diameter ${\rm diam}(\Omega)$ and inradius ${\rm inr}(\Omega)$, such latter quantity being defined as the supremum of radii of balls included in $\Omega$:
\[
\frac{\pi^2}{{\rm diam}(\Omega)^2}< \lambda_2-\lambda_1\le \frac{d \pi^2}{{\rm inr}(\Omega)^2}.
\]
This bounds can be somewhat improved when further geometrical properties of $\Omega$ hold, \cite{Smits}. Notice that, by taking $\Omega$ to be a rectangle of sides $L$ and $L^{-1}$ with $L$ large one explicitly computes $\lambda_2-\lambda_1=3\pi^2/L^2$, and $L$ is close to be the diameter of $\Omega$. A lower bound of the form $\lambda_2-\lambda_1>3\pi^2/{\rm diam}(\Omega)^2$ is conjectured to be the sharp one.

\subsection{Facts on the Elliptic Problem}

As mentioned above, the stabilization of the solutions $v\ge 0$ of the
transformed evolution problem (\ref{FDE.Problem.Domain.Fine.Log})
leads in a natural way to the consideration of the associated
stationary solutions, i.e., the solutions of the following
elliptic problem\label{Section.Elliptic}
\begin{equation*}\label{FDE.Elliptic.Problem1}
\begin{split}
\left\{\begin{array}{lll}
-\Delta (S^m)={\bf c}\,S & ~ {\rm in}~ \Omega,\\
S(x)>0& ~{\rm for}~ x\in\Omega,\\
S(x)=0 & ~{\rm for}~ x\in\partial\Omega,
\end{array}\right.
\end{split}
\end{equation*}
where $m_s<m<1$ and $\Omega\subset\RR^d$ is an open connected
domain with sufficiently smooth boundary. Using the new variable $V=S^m$ and putting
$p=1/m>1$  the latter problem can be written in the more popular semilinear elliptic form
\begin{equation*}\label{5.5}
-\Delta V={\bf c}V^p \quad \mbox{in } \ \Omega, \qquad V=0 \quad \mbox{on } \ \partial\Omega
\end{equation*}
Note that our restriction $m>m_s$ is the exact condition that makes the last problem subcritical, $p<p_s$.

\medskip

\noindent$\bullet~${\sl Existence of positive classical solutions}

\noindent The question of existence and regularity is well
understood in its basic features:

\noindent(a) if $0<m<1$ for $d\le 2$ or if
$\frac{d-2}{d+2}=m_s<m<1$ for $d\ge 3$, then there exist positive
classical solutions to equation (see e.g. \cite{BH} and references
quoted therein, and also \cite{DKV}).

\noindent(b) if $0<m\le m_s$ and $d\ge 3$ then there are cases in
which the positive classical solution exists (e.g. if $\Omega$ is
an annulus) and cases in which it does not  exist (e.g., if
$\Omega$ is star-shaped) (see e.g. \cite{BH} and references quoted
therein).

We observe that the geometry of the domain plays a role in the
question of existence, but only in the subcritical case (b), which
is not considered in this paper. Since we assume that $m>m_s$, the
existence of at least one positive classical solution is always
guaranteed.

\medskip

\noindent$\bullet~${\sl Uniqueness.} In the supercritical case
$m>m_s$ considered here, {\sl the geometry of $\Omega$ plays a
role in the uniqueness problem}. For example, if $d=1$ or if  $d\ge 2$ and $\Omega$ is a ball, then
the solution is unique, cf. \cite{AY}. While when $d\ge 2$ and $\Omega$ is an annulus, then the
solution is unique only in the class of positive radial solutions, cf. \cite{NN}. However, there are
cases in which the solution is not unique, cf. \cite{NN, BN}.

\medskip

\noindent$\bullet~${\sl Regularity and boundary behaviour.} We state now the main bounds for \eqref{5.5}, with explicit constants, for all $1\le p<p_s$. They will give us explicit bounds for the constants $\widetilde{k}_0(p), \widetilde{k}_1(p)$ appearing in Theorem \ref{Elliptic.THM}. We remark that we already know that $\widetilde{k}_i(p)\to 1$ as $p\to 1$, but we have no explicit bounds for them. While providing below such bounds for all $1\le p<p_c$, we notice that the resulting estimates will not satisfy the above limiting property.
The proofs follow by using the arguments that can be found for example in \cite{GT} for the local bounds, or in \cite{GNN, dFLN} for the boundary estimates, and they will be published separately in \cite{BGV-Elliptic}.
\begin{thm}[Local Upper Estimates]\label{thm.local.upper}
Let $\Omega\subset\RR^d$ be a bounded domain, and let ${\bf c}>0$. Let $V$ be a local weak (sub-)solution in $B_{R_0}\subset\Omega$ to $-\Delta V = {\bf c} V^p$, with $1\le p<p_s=2^*-1=(d+2)/(d-2)$. Then for any $R_\infty<R_0$  the following bound holds true:
\begin{equation}\label{upper.p+1}
\|V\|_{\infty,R_\infty}
    \le h_{1,p}\;\|V\|_{\overline{q},R_{0}}^{\frac{2\overline{q}}{2\overline{q}-d(p-1)}}\qquad\mbox{for any}\qquad \frac{d(p-1)}{2} < \overline{q}
\end{equation}
where the constant $h_{1,p}$ depends on $d,\,\overline{q},\,p,\,\mathcal{S}_2,\,R_0, R_\infty$ and can be explicitly calculated as  in {\rm \cite{BGV-Elliptic}}.
\end{thm}
\begin{thm}[Local Lower Estimates]\label{thm.local.lower}
Let $\Omega\subset\RR^d$ be a bounded domain, and let ${\bf c}>0$. Let $V$ be a local weak solution in $B_{R_0}\subset\Omega$ to $-\Delta V = {\bf c} V^p$, with $1\le p<p_s=2^*-1=(d+2)/(d-2)$. Then for any $\varepsilon>0$ and for any
\[
0<q\le\frac{2^{\frac{d-2}{2}}}{d\omega_d^2[\ee(d-1)+\varepsilon]}
\]
the following bound holds true
\begin{equation}\label{lower.local.thm}
\inf_{x\in B_{R_\infty}}V(x) =\|V\|_{-\infty,R_\infty}
    \ge h_0\frac{\|V\|_{q,R_{0}}}{|B_{R_0}|^{\frac{1}{q}}}\,.
\end{equation}
where the constant $h_0$ depends on $d,\,q,\,\varepsilon,\,\mathcal{S}_2,\,R_0, R_\infty$ and can be explicitly calculated as in {\rm \cite{BGV-Elliptic}}.
\end{thm}
By means of these upper and lower bounds one can prove quantitative Harnack estimates.
\begin{thm}[Harnack inequality for $1\le p<p_c$]\label{Harnack.ell.pc}
Let $\Omega\subset\RR^d$ be a bounded domain, and let ${\bf c}>0$\,. Let $V$ be a local weak solution in $B_{R_0}\subset\Omega$ to $-\Delta V = {\bf c} V^p$, with $1\le p < p_c=d/(d-2)$, and assume that $\|V\|_{p+1,R_0}\le K_p$. Then the following bound holds true for all $R_\infty<R_0$:
\begin{equation}\label{harnack.elliptic.pc}
\sup_{B_{R_\infty}}V(x)\le \mathcal{H}_p\inf_{B_{R_\infty}}V(x)
\end{equation}
where the constant $\mathcal{H}_p$ depends on $d,\,p,\,\mathcal{S}_2,\,R_0, R_\infty, K_p$ and can be explicitly calculated as in {\rm \cite{BGV-Elliptic}}.
\end{thm}
We now compare solutions corresponding to different $p$ and ${\bf c}$, and this can be done for $1\le p<p_c$\,, since we need the quantitative Harnack inequalities of Theorem \ref{Harnack.ell.pc}, that hold only in that range of $p$. We recall that we are now choosing ${\bf c}=1/[(1-m)T]=p/[(p-1)T]$ so that by Proposition \ref{prop.C} we have that ${\bf c}\to\lambda_1$ as $p\to 1$. Hence $\mathcal{H}_p$ in the above Theorem has a finite limit $\mathcal{H}_1$ as $p\to 1$.
\begin{thm}\label{prop.k0.k1.Appendix}
Let $U_p$ be a weak solution to the Elliptic problem {\rm \ref{ell.prob.p}}, without any assumption on $\lambda_p>0$, and let $1\le p<p_c=d/(d-2)$. There exists a $\delta(\Omega)=\delta>0$ independent of $U_p$, such that
\begin{equation}
\underline{k}_0(p) \le \frac{U_p(x)}{\Phi_1(x)}\le \overline{k}_1(p),
\end{equation}
where
\begin{equation}
\underline{k}_0(p)=\frac{1}{2\mathcal{H}_p}\frac{h_{1,p}}{h_{1,1}} \left[1-\frac{\delta}{R_1+\delta}\right]^{d-2},
\qquad
\overline{k}_1(p):= 2\mathcal{H}_1\frac{h_{1,p}}{h_{1,1}}\left[1-\frac{\delta}{R_1+\delta}\right]^{-(d-2)}
\end{equation}
and $\overline{c}_p$ are the constants in the upper bounds of Theorem {\rm \ref{thm.local.upper}} and $\mathcal{H}_p$ is the constant in the Harnack inequality of Theorem {\rm\ref{Harnack.ell.pc}}.
\end{thm}

\subsection{Maximum and comparison principles on small sets}
The maximum and comparison principle do not hold in general for solutions to nonlinear elliptic equations. This is an important characteristic of elliptic equation in general and does not necessarily depend on the nonlinearity. Indeed in the linear case, if one consider the Dirichlet problem for the equation $-\Delta u=\lambda u$ with $\lambda>\lambda_1$: it happens for instance that for $\lambda=\lambda_2>\lambda_1>0$ the corresponding second eigenfunction $\Phi_2$ has at least a change of sign, hence no global maximum nor comparison principle is allowed to hold.

In any case, we can still prove a (local) maximum and comparison principle on small sets: we are going to extend  to our framework an idea originally due to Serrin, see for example the book \cite{Serrin} where this idea is applied here to a different class of nonlinear elliptic equations. We just state the Theorem here, a complete proof will appear separately in \cite{BGV-Elliptic}.
\begin{thm}[Comparison with supersolutions on small sets]\label{thm.upper.comp}
Let $B\subset \RR^d$ be a bounded connected domain, let $p\ge 1$, $\lambda>0$ and
\[\left\{
\begin{array}{lll}
-\Delta u=\lambda\,u^p&\mbox{in }B\\
-\Delta \overline{u}\ge \lambda\,\overline{u}^p&\mbox{in }B\\
\overline{u}\ge u & \mbox{on }\partial B\\
0\le u,\overline{u}\le M &\mbox{in }\overline{B}\\
\end{array}\right.
\]
and assume that $|B|<\omega_d/\left(2p\,\lambda\,M^{p-1}\right)^d$. Then, we have that $\overline{u}\ge u$ in $\overline{B}$.
\end{thm}

\medskip

\noindent {\large \sc Acknowledgment}

\noindent  The first and third author have been partially funded by Project MTM2008-06326-C02-01 (Spain) and European Science Foundation Programme ``Global And Geometric Aspects of
Nonlinear Partial Differential Equations''. All authors acknowledge a contribution by the 2008 Spain-Italy research initiative HI2008-0178.

\small
\bibliographystyle{amsplain} 

\end{document}